\documentclass[a4paper]{amsart}

\usepackage[utf8]{inputenc}
\usepackage[british]{babel}
\usepackage{amsfonts,amsmath,amssymb,amsthm,mathtools,esint}
\usepackage{tikz}
\usepackage{enumitem}
\usepackage[hidelinks]{hyperref}
\usepackage{cleveref}
\crefname{subsection}{subsection}{subsections}
\numberwithin{equation}{section}
\usepackage{xcolor}
\definecolor{dark-gray}{gray}{0.3}
\usepackage{caption}

\usepackage[backend=bibtex, giveninits=true, isbn=false, doi=false, url=false, maxnames = 10]{biblatex}
\DeclareFieldFormat[article]{title}{{#1}} %
\DeclareFieldFormat[book]{title}{{#1}}
\DeclareFieldFormat[inproceedings]{title}{{#1}}
\DeclareFieldFormat[incollection]{title}{{#1}}
\DeclareFieldFormat{pages}{{#1}} %
\AtBeginBibliography{\footnotesize}
\renewbibmacro{in:}{%
	\ifentrytype{article}
	{}
	{\bibstring{in}%
		\printunit{\intitlepunct}}}

\newcommand{\diam}{\mathrm{diam}}
\newcommand{\R}{\mathbb{R}}

\renewcommand{\d}[1]{\,\mathrm{d}#1}
\renewcommand{\div}{\mathrm{div}}
\newcommand{\I}{\mathrm{I}}
\newcommand{\vertiii}[1]
{{\left\vert\kern-0.25ex\left\vert\kern-0.25ex\left\vert #1 
\right\vert\kern-0.25ex\right\vert\kern-0.25ex\right\vert}}
\newcommand{\Pkab}[2]{P_{#1}^{#2, 0}}

\newcommand{\Fcal}{\mathcal{F}}
\newcommand{\Tcal}{\mathcal{T}}

\newcounter{cnt}
\newcommand{\newcnst}{%
	\refstepcounter{cnt}%
	\ensuremath{C_{\thecnt}}}
\newcommand{\cnst}[1]{\ensuremath{C_{\ref{#1}}}}

\newcommand{\GrRec}{\mathcal{G}}
\newcommand{\PotRec}{\mathcal{R}}
\newcommand{\RT}{\mathrm{RT}}
\newcommand{\pw}{\mathrm{pw}}

\newtheorem{theorem}{Theorem}[section]
\newtheorem{lemma}[theorem]{Lemma}

\newtheorem{corollary}[theorem]{Corollary}
\newtheorem*{conjecture}{Conjecture}
\theoremstyle{remark}
\newtheorem{remark}[theorem]{Remark}
\newtheorem{example}[theorem]{Example}

\begin{document}

\title[HHO for lower eigenvalue bounds]{Adaptive hybrid high-order method for guaranteed lower eigenvalue bounds}

\author[C.~Carstensen, B.~Gr\"a\ss le, N.~T.~Tran]{Carsten Carstensen \and Benedikt Gr\"a\ss le \and Ngoc Tien Tran}

\thanks{This work has been supported by the \emph{German Research Foundation} (DFG) in the Priority Program 1748 \emph{Reliable simulation techniques in solid mechanics: Development of
non-standard discretization methods, mechanical and mathematical analysis} CA 151/22-2 and 
under Germany's Excellence Strategy – The Berlin Mathematics Research Center MATH+ (EXC-2046/1, project ID: 390685689) as well as
by the European Union's Horizon 2020 research and innovation programme (project DAFNE grant agreement No.~891734 and project RandomMultiScales grant agreement No.~865751).}

\address[C.~Carstensen, B.~Gr\"a\ss le]{%
	Institut f\"ur Mathematik,
	Humboldt-Universit\"at zu Berlin,
	10117 Berlin, Germany}
\email{cc@math.hu-berlin.de, graesslb@math.hu-berlin.de}
\address[N.~T.~Tran]{%
	Institut f\"ur Mathematik,
	Universit\"at Augsburg,
	86159 Augsburg, Germany}
\email{ngoc1.tran@uni-a.de}



\begin{abstract}
	The higher-order guaranteed lower eigenvalue bounds of the Laplacian 
	in the recent work by Carstensen, Ern, and Puttkammer [Numer. Math. 149, 2021]
	require a parameter $C_{\mathrm{st},1}$ that is found {\em not} robust as the polynomial 
	degree $p$ increases.
	This is related to the $H^1$ stability bound of the $L^2$ projection onto polynomials of degree at most $p$ and its
	growth $C_{\rm st, 1}\propto (p+1)^{1/2}$ as $p \to \infty$. A similar estimate for the Galerkin projection holds with a $p$-robust constant $C_{\mathrm{st},2}$ and $C_{\mathrm{st},2} \le 2$ for right-isosceles triangles.
	This paper utilizes the new inequality with the constant $C_{\mathrm{st},2}$
	to design a modified hybrid high-order (HHO) eigensolver that 
	directly computes guaranteed lower eigenvalue bounds under the idealized 
	hypothesis of exact solve of the generalized algebraic eigenvalue problem
	and a mild explicit condition on the maximal mesh-size in the simplicial mesh. 
	A key advance is a $p$-robust parameter selection.
	
	The analysis of the new method with a different fine-tuned volume stabilization allows for a~priori quasi-best approximation and improved $L^2$ error estimates as well as a stabilization-free 
	reliable and efficient a posteriori error control. The associated adaptive mesh-refining algorithm performs superior
	in computer benchmarks with striking numerical evidence for optimal higher empirical convergence rates.
\end{abstract}

\maketitle

\section{Introduction}
This paper proposes and analyzes a new hybrid high-order (HHO) eigensolver for the direct computation of guaranteed lower eigenvalue bounds (GLB) for the Laplacian.

\subsection{Three categories of GLB}
The min-max principle enables guaranteed upper eigenvalue bounds but prevents a direct computation of a GLB by a conforming approximation in a Rayleigh quotient.
So GLB shall be based on nonconforming finite element methods (FEM),
on modified mass and/or stiffness matrices (with reduced integration or fine-tuned stabilization terms), or on further post-processing.
The last decade has seen a few GLB we group in three categories
(i)--(iii).

\begin{enumerate}[wide]
	\item[(i)] The a posteriori error analysis for symmetric second-order elliptic eigenvalue problems started with \cite{Verfuerth1986,Larson2000,DuranPadraRodriguez2003} under the (unverified) hypothesis of a sufficiently small mesh-size. With additional a priori information  
	on spectral gaps, the latest a posteriori post-processings \cite{CancesDussonMadayVohralik2017,CancesDussonMadayVohralik2018,CancesDussonMadayVohralik2020} provide GLB.
	
	\item[(ii)] Classical nonconforming FEM \cite{CarstensenGallistl2014,CarstensenGedicke2014,Liu2015}
	and mixed FEM \cite{Gallistl2023}
	allow for the
	GLB $\lambda_h/(1+\delta\lambda_h)\le \lambda$ with the discrete eigenvalue $\lambda_h$ and a known parameter
	$\delta\propto h_{\max}^2$ in terms of the maximal mesh-size $h_{\max}$.
	On the positive side, the GLB provides unconditional information on the exact eigenvalue $\lambda$ from the computed
	discrete eigenvalue $\lambda_h$.
	On the negative side, the global parameter $h_{\max}$ can spoil a very accurate  approximation $\lambda_h$ in this
	GLB and is
 of lowest-order only.
	A fine-tuned stabilization of the classical
	nonconforming FEM in \cite{CP21b}, however, provides a first (but low-order) remedy
	of the third category.
	\item[(iii)] 
		Higher-order hybrid discontinuous Galerkin (HDG) or HHO discretizations \cite{CarstensenZhaiZhang2020,CarstensenErnPuttkammer2021} can compute direct GLB $\lambda_h\le \lambda$ under the sufficient condition	(e.g., in \cite{CarstensenErnPuttkammer2021} for the HHO method and the Laplacian) 
	\begin{align}
		\sigma_1^2\beta +\kappa^2 h_{\max}^2\min\{\lambda,\lambda_h\}\le \alpha
		\label{eq:cond}
	\end{align}
	with (universal or computed) constants $\sigma_1, \kappa$ and known parameters $0<\alpha<1, 0<\beta<\infty$ (selected in the discrete scheme).
	If the exact Dirichlet eigenvalue $\lambda$ of number $j \in \mathbb{N}_0$ of the Laplace operator and the corresponding discrete eigenvalue $\lambda_h$ satisfy \eqref{eq:cond}, then $\lambda_h \le \lambda$ is a GLB.
	The two-fold use of \eqref{eq:cond} is a priori or a posteriori.
	First, given an upper bound $\mu \ge \lambda > 0$ of $\lambda$ (e.g., by some conforming approximation or post-processing), \eqref{eq:cond} provides  
	an upper bound $h_{\max}^2 \le (\alpha - \sigma_1^2\beta)/(\kappa^2\mu)$ for the maximal initial mesh-size. This condition is sufficient for \eqref{eq:cond} and guarantees a priori that $\lambda_h \le \lambda $. Second,
	\eqref{eq:cond} may be checked a posteriori for any computed value $\lambda_h$. Then
	$\sigma_1^2\beta +\kappa^2 h_{\max}^2 \lambda_h \le \alpha$ implies \eqref{eq:cond} and so, $\lambda_h \le \lambda$. 
\end{enumerate}
This paper presents a new HHO eigensolver of the third category.
\begin{figure}
	\centering
	\includegraphics[]{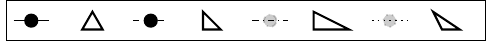}
	\includegraphics[]{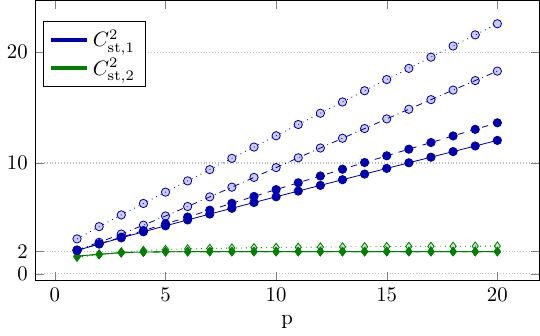}
	\caption{Approximations %
		of $C_{\rm st,1}^2$ and $C_{\rm st,2}^2$ as a function of the polynomial degree $p$ on the equilateral triangle
		\protect\includegraphics[height=.8em]{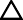}, the right-isosceles triangle \protect\includegraphics[height=.8em]{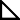},
		\protect\includegraphics[height=.8em]{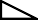}$\coloneqq \mathrm{conv}\{(0,0), (1.5,0), (0,1)\}$, and
		\protect\includegraphics[height=.8em]{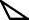}$\coloneqq\mathrm\{(0,0), (1,0), (-1/2, 1)\}$.}
	\label{fig:C_p}
\end{figure}

\subsection{Motivation and outline of \Cref{sec:stability_estimates}}
The constants $\sigma_1^2 \coloneqq C_\textrm{st,1}^2 - 1$ and $\kappa \coloneqq C_P C_\textrm{st,1}$ in \eqref{eq:cond} depend on the Poincar\'e constant $C_P \leq 1/\pi$ and a stability constant $C_\textrm{st,1}$. The latter has to be contrasted with the
constant $C_{\rm st, 2}$, where $C_{\rm st, 1}$ and $C_{\rm st, 2}$ are the best possible constants in the stability estimates
\begin{align}
	\|\nabla (1 - \Pi_{p+1}) f\|_{L^2(T)} &\leq C_\textrm{st,1}\|(1 - \Pi_p)\nabla f\|_{L^2(T)}\quad\text{for all }f\in H^1(T),
	\label{ineq:stability-estimate-1}\\
	\|\nabla (1 - G_{p+1}) f\|_{L^2(T)} &\leq C_\textrm{st,2}\|(1 - \Pi_p)\nabla f\|_{L^2(T)}\quad\text{for all }f\in H^1(T)
	\label{ineq:stability-1}
\end{align}
in a given simplex $T \subset \R^n$ with
the (component-wise) $L^2$ projection $\Pi_m$ and the Galerkin projection $G_m$ onto polynomials of total degree at
most $m\in\mathbb N_0$.
The two constants $C_\mathrm{st,1}$ and $C_\mathrm{st,2}$ are independent of the diameter $h_T \coloneqq \mathrm{diam}(T)$ of $T$, but might depend on the shape of $T$ and the polynomial degree $p$.
\Cref{fig:C_p} illustrates the behaviour of $C_{\rm st, 1}$ and $C_{\rm st, 2}$ for different triangular shapes and various polynomial degrees $p$.
Section \ref{sec:stability_estimates} investigates the $p$-robustness of $C_{\rm st, 2}$ and reveals that $C_{\rm st, 2}\leq C_{\rm st,1}\propto\sqrt{p+1}$ tends to
infinity as $p\to \infty$, while we conjecture $C_{\rm st,2}\leq \sqrt 2$ for triangles $T$ with maximum interior angle
$\omega\leq \pi/2$.
Notice that a large constant $C_\mathrm{st,1}$ leads to a large $\sigma_1$ in \eqref{eq:cond} and so, $\alpha < 1$ enforces small $\beta$ and restricts the GLB to very fine meshes.
The main motivation 
of this work arises from the convenient bound $C_{\rm st, 2}\leq \sqrt2$: Can we design a discretization method of the
third category (iii) based on $\sigma_2 \coloneqq C_P C_\textrm{st,2}\leq\sqrt2/\pi$ in \eqref{eq:cond}?
\subsection{A modified HHO method and outline of \Cref{sec:HHO}}
This paper provides an affirmative answer with a new fine-tuned stabilization in a modified HHO scheme in
\Cref{sec:HHO} and a new criterion
\begin{align}
	\sigma_2^2 \max\{\beta, h_{\max}^2\min\{\lambda,\lambda_h\}\} \le \alpha
	\label{eq:cond2}
\end{align}
sufficient for the GLB $\lambda_h\le\lambda$.
One advantage of \eqref{eq:cond2} over \eqref{eq:cond} is the straight-forward and $p$-robust parameter selection $\beta \coloneqq
\alpha/\sigma_2^2$.
It turns out that $\sigma_2 \leq \kappa$ and so \eqref{eq:cond2}
improves on \eqref{eq:cond} in the sense that $\sigma_2^2 h_{\max}^2\lambda \leq \alpha$ holds on much coarser triangulations for higher polynomial degrees $p$.

Given a bounded polyhedral Lipschitz domain $\Omega \subset \R^n$, let
$V \coloneqq H^1_0(\Omega)$ denote the Sobolev space endowed with the energy scalar product
$a(u,v) \coloneqq (\nabla u, \nabla v)_{L^2(\Omega)}$ and the $L^2$ scalar product $b(u,v) \coloneqq (u,v)_{L^2(\Omega)}$ for all $u,v \in V$.
This paper considers the model problem that seeks an eigenpair $(\lambda, u) \in \R_+ \times (V\setminus\{0\})$ such that
\begin{align}
	a(u, v) = \lambda b(u,v) \text{ for any } v \in V.
	\label{pr:eigenvalue-problem}
\end{align}
The HHO methodology has been introduced in \cite{DiPietroErnLemaire2014,DiPietroErn2015} and is related to HDG and nonconforming virtual element methods \cite{CockburnDiPietroErn2016}.
Given a regular triangulation $\Tcal$ into simplices,
the ansatz space $V_h = P_{p+1}(\Tcal) \times P_p(\Fcal(\Omega))$ consists of piecewise polynomials of (total) degree at most
$p+1$ attached to the simplices and piecewise polynomials of degree at most $p$ attached to the interior faces.
Two reconstruction operators link the two components of $v_h \in V_h$:
The potential reconstruction $\PotRec v_h \in P_{p+1}(\Tcal)$ provides a discrete approximation to $v$ in the space of piecewise polynomials $P_{p+1}(\Tcal)$ of degree at most $p+1$.
The gradient reconstruction $\GrRec v_h \in \RT_p^\pw(\Tcal)$ approximates the gradient $\nabla v$ in the space of piecewise Raviart-Thomas functions $\RT_p^\pw(\Tcal)$ \cite{AbbasErnPignet2018,DiPietroDroniouManzini2018}.
Let $S v_h \coloneqq v_\Tcal - \PotRec v_h \in P_{p+1}(\Tcal)$ for any $v_h = (v_\Tcal, v_\Fcal) \in V_h$ denote the
additional cell-based stabilization.
Given positive parameters $0 < \alpha < 1$ and $0 < \beta < \infty$, the bilinear forms $a_h : V_h \times V_h \to \R$
and $b_h: V_h \times V_h \to \R$ read
\begin{align}
	a_h(u_h,v_h) &\coloneqq (\GrRec u_h, \GrRec v_h)_{L^2(\Omega)} - \alpha ((1-\Pi_p)\GrRec u_h, (1-\Pi_p)\GrRec v_h)_{L^2(\Omega)}\label{def:ah}\\
	&\qquad + \beta (h_\mathcal{T}^{-2} S u_h, S v_h)_{L^2(\Omega)}\nonumber,\\
	b_h(u_h,v_h) &\coloneqq (u_\Tcal, v_\Tcal)_{L^2(\Omega)} \text{ for any } u_h = (u_\Tcal,u_\Fcal), v_h = (v_\Tcal, v_\Fcal) \in V_h.\label{def:bh}
\end{align}
The discrete eigenvalue problem seeks $(\lambda_h,u_h)\in\mathbb{R}^+\times (V_h\setminus\{0\})$ with
\begin{align}
	a_h(u_h,v_h)=\lambda_h b_h(u_h,v_h) \text{ for all } v_h\in V_h. \label{eq:disEVP}
\end{align}
The definitions of $\PotRec$, $\GrRec$, and further details follow in \Cref{sec:HHO} below.

\subsection{GLB with $p$-robust parameters and outline of \Cref{sec:GLB}}%
	\label{sub:GLB with $p$-robust parameters and outline of Cref{sec:GLB}}
The discrete bilinear form $a_h$ from \cite{CarstensenErnPuttkammer2021} with parameter $C_{\rm st,
1}\propto\sqrt{p+1}$ utilizes the different stabilization $\beta(\nabla_\pw S u_h, \nabla_\pw S v_h)_{L^2(\Omega)}$ instead of $\beta(h_\Tcal^{-2} S u_h, S v_h)_{L^2(\Omega)}$ in \eqref{def:ah}.
The two stabilizations are locally equivalent, but the innovative difference is that the parameter selection in the new scheme circumvents an inverse
inequality, and rather builds it into the scheme.
\Cref{sec:GLB} verifies the sufficient condition \eqref{eq:cond2} for exact GLB under the assumption of exact solve.

\subsection{A priori error analysis of the new scheme and outline of \Cref{sec:a-priori}}%
\label{sub:Analysis of the new scheme}
A quasi-best approximation
for the source problem \cite{ErnZanotti2020} allows for quasi-best approximation results in \Cref{thm:L2-error-control}
for a simple eigenvalue $\lambda$, namely
\begin{align}
	|\lambda - \lambda_h| + a_h(\I u - u_h, \I u - u_h) + h_{\max}^{-2s}\|u - u_\Tcal\|_{L^2(\Omega)}^2&\nonumber\\
	\leq \cnst{cnst:best-approximation} \min_{v_{p+1} \in P_{p+1}(\Tcal)}\|\nabla_\pw(u - v_{p+1})\|_{L^2(\Omega)}^2&\label{ineq:quasi-optimal-eigenvalue}
\end{align}
with a positive constant $\newcnst\label{cnst:best-approximation}$ and the minimum 
$0 < s \leq 1$ of the index of elliptic regularity and one;
the HHO interpolation $\I : V \to V_h$ is recalled in \Cref{sec:HHO-methodology} below.
Compared to earlier results in \cite{CaloCicuttinDengErn2019,CarstensenErnPuttkammer2021},
\eqref{ineq:quasi-optimal-eigenvalue} provides an additional positive power $s$ of $h_{\max}$ in the $L^2$ error. This
is important as it eventually enables the absorption of higher-order terms in the a~posteriori error analysis.

\subsection{Stabilization-free a~posteriori error analysis and outline of \Cref{sec:a-posteriori}}
Let $p_h \coloneqq \Pi_p \GrRec u_h \in P_p(\Tcal;\R^n)$ denote the $L^2$ projection of the gradient reconstruction $\GrRec u_h \in \RT_p^\pw(\Tcal)$ onto the space of vector-valued piecewise polynomials $P_p(\Tcal;\R^n)$.
For any $T \in \Tcal$ of volume $|T|$, define
\begin{align}
	&\eta^2(T) \coloneqq |T|^{2/n} \big(\|\div\,p_h + \lambda_h u_\Tcal\|_{L^2(T)}^2 + \|\mathrm{curl}\,p_h\|^2_{L^2(T)}\big)
	\label{def:eta-local}\\
	&\quad  + |T|^{1/n}\Big(\sum_{F \in \Fcal(T) \cap \Fcal(\Omega)} \|[p_h \cdot \nu_F]_F\|_{L^2(F)}^2 + \sum_{F \in \Fcal(T)} \|[p_h \times \nu_F]_F\|_{L^2(F)}^2\Big)
	\nonumber
\end{align}
with the normal jump $[p_h \cdot \nu_F]_F$ and the tangential jump $[p_h \times \nu_F]_F$ of $p_h$ across a side $F$ of
$T$. \Cref{thm:a-posteriori} asserts reliability and efficiency of the error estimator
$\eta^2 \coloneqq \sum_{T \in \Tcal} \eta^2(T)$ 
for sufficiently small mesh-sizes $h_{\max}$ in the sense that
\begin{align}
	C_{\mathrm{eff}}^{-1} \eta \leq |\lambda - \lambda_h| + a_h(\I u - u_h, \I u - u_h) + \|\nabla
	u - p_h\|^2_{L^2(\Omega)} \leq C_{\mathrm{rel}} \eta.
	\label{ineq:a-posteriori}
\end{align}

\subsection{Adaptive mesh-refining algorithm and outline of \Cref{sec:numerical-examples}}
Three 2D computer experiments in \Cref{sec:numerical-examples} provide striking numerical evidence that the criterion \eqref{eq:cond2} indeed leads to confirmed lower eigenvalue bounds.
The adaptive mesh-refining algorithm driven by the refinement indicator $\eta$ from \eqref{def:eta-local} recovers the optimal convergence rates of the eigenvalue error $\lambda - \lambda_h$ in all numerical benchmarks with singular eigenfunctions.
This is the first time that $p$-robust higher-order GLB of the third category are displayed.

\subsection{General notation}
Standard notation for Lebesgue and Sobolev function spaces applies throughout this paper. In particular,
$(\bullet,\bullet)_{L^2(\omega)}$ denotes the $L^2$ scalar product and $H(\div,\omega)$ is the space of Sobolev
functions with weak divergence in $L^2(\omega)$ for a domain $\omega\subset \R^n$.
Recall the abbreviation $V \coloneqq H^1_0(\Omega)$ for the space of Sobolev functions,
endowed with the energy scalar product
$a(u,v) \coloneqq (\nabla u, \nabla v)_{L^2(\Omega)}$ and the $L^2$ scalar product $b(u,v) \coloneqq (u,v)_{L^2(\Omega)}$ for all $u,v \in V$.

For a subset $M \subset \R^n$ of diameter $h_M$, let $P_p(M)$ denote the space of polynomials of maximal (total) degree $p$ regarded as functions defined in $M$.
Given a simplex $T \subset \R^n$, the space of Raviart-Thomas finite element functions reads
\begin{align*}
	\RT_p(T) &\coloneqq P_p(T;\R^n) + x P_p(T) \subset P_{p+1}(T;\R^n).
\end{align*}
The Galerkin projection $G \coloneqq G_{p+1}: H^1(T) \to P_{p+1}(T)$ maps $f \in H^1(T)$ to the unique solution $G f$ to $\Pi_0 G f = \Pi_0 f$ and
\begin{align}
	(\nabla G f, \nabla p_{p+1})_{L^2(T)} = (\nabla f,\nabla p_{p+1})_{L^2(T)} \text{ for all } p_{p+1} \in P_{p+1}(T)
	\label{def:Galerkin}
\end{align}
with the convention $H^1(T) \coloneqq H^1(\mathrm{int}(T))$ for the interior $\mathrm{int}(T) = T^\circ$ of $T$.
The Poincar\'e constant $C_P$ bounds $\|(1-\Pi_0)f\|_{L^2(T)}\leq C_Ph_T\|\nabla
f\|_{L^2(T)}$ for all $f\in H^1(T)$.
In 2D, $C_P \leq 1/j_{11} = 0.260980$ with the first positive root of the Bessel function $J_1$
\cite{LaugesenSiudeja2010} and $C_P \leq 1/\pi$ in any space dimension \cite{PayneWeinberger1960,Bebendorf2003}.
The context-sensitive notation $|\bullet|$ may denote the absolute value of a scalar, the Euclidean norm of a vector, the length of a side, or the volume of a simplex.
The notation $A \lesssim B$ abbreviates $A \leq CB$ for a generic constant $C$ independent of the mesh-size and $A \approx B$ abbreviates $A \lesssim B \lesssim A$. 
Throughout this paper, $C_1, \dots, \cnst{cnst:a-posteriori-absorption-small-mesh-size}$ denote positive constants independent of the mesh-size.

\section{Stability estimates}\label{sec:stability_estimates}
This section discusses the behaviour of the constants $C_{\rm st, 1}, C_{\rm st, 2}$ from
\eqref{ineq:stability-estimate-1}--\eqref{ineq:stability-1} as $p\to\infty$
and the computation of %
$\sigma_2\coloneqq C_PC_{\rm st, 2}$ with the Poincar\'e constant $C_P$ in \eqref{eq:cond2} that arises from the stability estimates in
\Cref{lem:stability_estimates} below.

\subsection{Stability constants and estimates}%
\label{sub:Stability constants and estimates}
The following theorem asserts that $C_{\rm st, 2}$ is $p$-robust (and small in general, see Figure \ref{fig:C_p}) whereas $C_{\rm st, 1}\to
\infty$ as $p\to \infty$.

\begin{theorem}
	For any simplex $T$, there exist positive constants $1 \leq C_{\rm st, 2}\leq C_{\rm st, 1}$ independent of the diameter
	$h_T$ such that \eqref{ineq:stability-estimate-1}--\eqref{ineq:stability-1} hold.
In $n=2,3$ space dimensions, %
$
	C_{\rm st, 1}\approx \sqrt{p+1}$ and $%
	C_{\rm st, 2}\approx 1  $ as $p\to\infty$.
\end{theorem}
\begin{proof}
	The existence of the constants $1 \leq C_{\rm st, 1} \leq C_{\rm st, 2}$ %
	follows from \cite[Theorem 3.1]{CarstensenZhaiZhang2020}; cf.~Appendix A for further details.
	The technical proof of the $p$-robustness of $C_\mathrm{st,2}$ involves a linear bounded operator $R^{\mathrm{curl}} :
	H^{-1}(T;\R^{3}) \to L^2(T;\R^3)$ from \cite{Hiptmair2009,CostabelMcIntosh2010,MelenkRojik2020} and is carried out in Appendix B.
	The robustness holds for $n = 2$ with a simpler and hence omitted proof.
	The remaining parts of this proof concern the growth of $C_{\rm st, 1}$.
	Let $X\coloneqq H^1(T)/\R$ denote the Hilbert space with inner product $(\nabla\bullet,
	\nabla\bullet)_{L^2(T)}$ and note that
	$\mathrm{ker}(\nabla(1-\Pi_{p+1}))=\mathrm{ker}((1-\Pi_p)\nabla)=P_{p+1}(T)$.
	Since $\|(1-\Pi_p)\nabla\phi\|_{L^2(T)}\leq\vertiii{\phi}$ for every $\phi\in X$, the definition of the
	operator norm of the oblique projection $1-\Pi_{p+1}\in L(X;X)$ provides
	\begin{align*}
		\|1-\Pi_{p+1}\|\coloneqq \sup_{\phi\in X\setminus\{0\}}\frac{\vertiii{(1-\Pi_{p+1})\phi}}{\vertiii{\phi}}\leq \sup_{\phi\in
		X\setminus P_{p+1}(T)}\frac{\vertiii{(1-\Pi_{p+1})\phi}}{\|(1-\Pi_p)\nabla\phi\|_{L^2(T)}}=C_{\rm st, 1}.
	\end{align*}
	Kato's oblique projection lemma \cite{szyld_many_2006} for the Hilbert space $X$ leads to
	$\|\Pi_{p+1}\|=\|1-\Pi_{p+1}\| \,\leq\, C_{\rm st, 1}$
	and $(1-\Pi_{p+1})G= 0$ in $X$ for the Galerkin projection $G$ shows
	\begin{align*}
		\vertiii{(1-\Pi_{p+1})f}
		=\vertiii{(1-\Pi_{p+1})(1-G)f}
		\leq \|\Pi_{p+1}\|\;\vertiii{
		(1-G)f}&&\text{for any }f\in X.
	\end{align*}
	Since $\vertiii{
		(1-G)f}\leq C_{\rm st, 2} \|(1-\Pi_p)\nabla f\|_{L^2(T)}$ from \eqref{ineq:stability-1}, this
	proves $\|\Pi_{p+1}\|\leq C_{\rm st, 1}\leq
	C_{\rm st, 2}\|\Pi_{p+1}\|$.
	The growth $\|\Pi_{p+1}\|\approx \sqrt{p+1}$ is known for tensor-product domains and also holds for
	simplices in $n=2,3$ dimensions; see
	\cite{wurzer2010stability} and 
	\cite[Sec.\ 5]{melenk_stability_2013} for $\|\Pi_{p+1}\|\lesssim \sqrt {p+1}$
	and Appendix C for the proof of $\sqrt {p+1}\lesssim \|\Pi_{p+1}\|$.
\end{proof}
The Poincar\'e inequality with the Poincar\'e constant $C_P$ and \eqref{ineq:stability-1}
with $C_{\rm st, 2}\approx 1$ lead to a $p$-robust stability
estimate with $\sigma_2\coloneqq C_PC_{\rm st, 2}$.
\begin{lemma}[$p$-robust stability]\label{lem:stability_estimates}
	Any $f\in H^1(T)$, $T$ a simplex, and $p\in \mathbb{N}_0$ satisfy 
	\begin{align}
		\|h_T^{-1}(1 - G) f\|_{L^2(T)} &\leq \sigma_2\|(1 - \Pi_p)\nabla f\|_{L^2(T)}.\qed
		\label{ineq:stability-2}
	\end{align}
\end{lemma}
\subsection{Numerical comparison and conjecture}%
\label{sub:numerical comparison}

The following theorem considers the computation of guaranteed upper bounds of $C_{\rm st, 2}$ in $n=2,3 $ space
dimensions for a control of $\sigma_2$ in \eqref{ineq:stability-2}.
\begin{table}[h]
	\centering
	\begin{tabular}{|c|c|c|c|}
		\hline
		$p$ & $C_\mathrm{st,2}^2$ & $p$ & $C_\mathrm{st,2}^2$\\\hline
		$1$ & 1.59707221 & $6$ & 1.99368122\\
		$2$ & \multicolumn{1}{|l|}{1.75} & $7$ & 1.99787853\\
		$3$ & 1.91060394 & $8$ & 1.99911016\\
		$4$ & 1.95679115 & $9$ & 1.99969758\\
		$5$ & 1.98559893 & $10$ & 1.99987656\\\hline
	\end{tabular}
	\caption{The constant $C_\mathrm{st,2}=m_p$ on right-isosceles triangles.}
	\label{tab:C_k}
\end{table}
\noindent
Given $v \in H^1(T;\R^n)$ and $w \in H^1(T;\R^{2n-3})$, let $\mathrm{curl}\,v \coloneqq \partial_1 v_2 - \partial_2 v_1$
and $\mathrm{Curl}\,w \coloneqq (\partial_2 w, -\partial_1 w)^t$ for $n = 2$ and $\mathrm{curl}\,v \coloneqq (\partial_2
v_3 - \partial_3 v_2, \partial_3 v_1 - \partial_1 v_3, \partial_1 v_2 - \partial_2 v_1)^t$ and $\mathrm{Curl}\,w
\coloneqq \mathrm{curl}\,w$ for $n = 3$.
For any $g \in H^{-1}(T;\R^{2n-3})$ in the dual space of $H^1_0(T;\R^{2n-3})$ endowed with the operator norm
$\vertiii{\bullet}_{*}$,
let  $(-\Delta)^{-1} g \in H^1_0(T;\R^{2n-3})$ denote the weak solution to $-\Delta v = g$ in $T$ componentwise with $\vertiii{g}_* = \vertiii{(-\Delta)^{-1} g}$.

\begin{figure}
	\centering
	\includegraphics[]{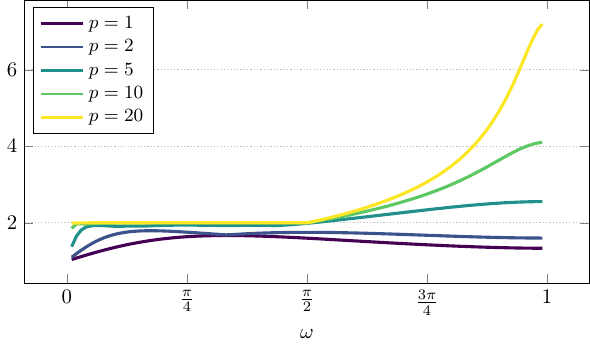}
	\caption{Dependence of $C_{\rm st,2}^2$ on the interior angle $\omega$ of the isosceles triangle
	$T=\mathrm{conv}\{(0,0), (1,0), (\cos(\omega),\sin(\omega))\}$}
	\label{fig:C_p_angle}
\end{figure}

The gradients $\nabla P_{p+1}(T)$ of polynomials $P_{p+1}(T)$ of degree at most $p+1$ form a subspace of $P_p(T;\R^n)$ and give rise to the $L^2$ orthogonal decomposition $P_p(T;\R^n) = Q_p \oplus \nabla P_{p+1}(T)$ with $Q_p \perp \nabla P_{p+1}(T)$ in $L^2(T;\R^n)$. Let $\mathrm{P}: P_p(T;\R^n) \to \nabla P_{p+1}(T)$ denote the $L^2$ orthogonal projection onto $\nabla P_{p+1}(T) \subset P_p(T;\R^n)$.
The bilinear forms $\mathfrak{a}: Q_p \times Q_p \to \R$ and $\mathfrak{b}: Q_p \times Q_p \to \R$ are defined, for any $q_p,r_p \in Q_p$, by
\begin{align}
	\mathfrak{a}(q_p,r_p) \coloneqq (q_p, r_p)_{L^2(T)} \quad\text{and}\quad
	\mathfrak{b}(q_p,r_p) \coloneqq ((-\Delta)^{-1} \mathrm{curl}\,q_p, \mathrm{curl}\,r_p)_{L^2(T)}.
	\label{def:a-b}
\end{align}
\begin{theorem}[stability constant]\label{thm:stability-constant}
	The maximal eigenvalue
	\begin{align}
		m_p^2 \coloneqq \max_{\substack{q_p \in P_p(T;\R^n)\\\mathrm{curl}\,q_p \neq 0}} \min_{v_{p+1} \in P_{p+1}(T)} \|q_p - \nabla v_{p+1}\|_{L^2(T)}^2/\vertiii{\mathrm{curl}\, q_p}_{*}^2
		\label{def:mk}
	\end{align}
 	of the eigenvalue problem
	\begin{align}
		\mathfrak{a}(q_p,r_p) = \lambda \mathfrak{b}(q_p,r_p) \quad\text{for all } r_p \in Q_p
		\label{def:mk-eigenvalue-problem}
	\end{align}
	leads to the upper bound $C_\mathrm{st,2} \leq \max\{1,m_p C_n\}$ for $\newcnst = 1$ and $\newcnst = 2/\sqrt 3$.
\end{theorem}
Notice that \eqref{def:mk-eigenvalue-problem} is a finite-dimensional eigenvalue problem and $(-\Delta)^{-1}q_p$ in $\mathfrak{b}(q_p,r_p)$ can be approximated by, e.g., a conforming FEM.
Numerical experiments below even suggest that the bound $C_{\rm st, 2}=m_p$ is exact in $n=2$ dimensions.
\begin{proof}
	If $p = 0$, $\nabla P_1(T) = P_0(T;\R^n)$ implies $\nabla G f = \Pi_0 \nabla f$ for all $f \in H^1(T)$, whence $C_\mathrm{st,2} = 1$.
	The remaining parts of the proof therefore assume $p \geq 1$.
	Given $f \in H^1(T)$, assume without loss of generality that $\nabla f \perp \nabla P_{p+1}(T)$ in $L^2(T;\R^n)$ (otherwise substitute $g \coloneqq f - Gf$ and observe that $\|(1 - \Pi_p) \nabla f\|_{L^2(T)} = \|(1 - \Pi_p)\nabla g\|_{L^2(T)}$).
	Throughout this proof, abbreviate $q_p \coloneqq \Pi_p \nabla f \in P_p(T;\R^n)$.
	A Helmholtz decomposition leads to $q_p = \nabla a + \mathrm{Curl}\,b$ with $a \in H^1(T)$ and $b \in H^1_0(T;\R^{2n-3})$.
	For any $v \in H^1_0(T;\R^{2n-3})$, the $L^2$ orthogonality $\mathrm{Curl}\,v \perp \nabla a$ in $L^2(T;\R^n)$, an integration by parts, and a Cauchy inequality prove
	\begin{align}
		\int_T v\cdot\mathrm{curl}\,q_p \d{x} &= \int_T q_p \cdot \mathrm{Curl}\,v \d{x}\nonumber\\
		&= \int_T \mathrm{Curl}\,b \cdot \mathrm{Curl}\,v \d{x} \leq \|\mathrm{Curl}\,b\|_{L^2(T)} \|\mathrm{Curl}\,v\|_{L^2(T)}.
		\label{ineq:proof-constant-integration-by-parts}
	\end{align}
	In 2D, $\|\mathrm{Curl}\,v\|_{L^2(\Omega)} = \vertiii{v}$ and in 3D, $\|\mathrm{Curl}\,v\|_{L^2(\Omega)} \leq 2
	\vertiii{v}/\sqrt{3}$. (The proof solely involves elementary algebra and is therefore omitted.) Hence,
	\eqref{ineq:proof-constant-integration-by-parts} implies
	\begin{align}
		\vertiii{\mathrm{curl}\,q_p}_* = \sup_{v \in H^1_0(T;\R^{2n-3})\setminus\{0\}} \int_T v\cdot\mathrm{curl}\,q_p \d{x}/\vertiii{v} \leq C_n\|\mathrm{Curl}\,b\|_{L^2(T)}.
		\label{ineq:proof-constant-curl-dual-norm}
	\end{align}
	(Notice that $\vertiii{q_p}_* = \|\mathrm{Curl}\,b\|_{L^2(T)} = \vertiii{b}$ in 2D). 
	Since $\nabla P_{p+1}(T) \subset P_p(T;\R^n)$, the best approximation of $q_p$ in $\nabla P_{p+1}(T)$ satisfies the $L^2$ orthogonality $q_p \perp \nabla P_{p+1}(T)$. This and the Pythagoras theorem provide
	\begin{align*}
		\min_{v_{p+1} \in P_{p+1}(T)} \|q_p - \nabla v_{p+1}\|_{L^2(T)}^2 = \|q_p\|_{L^2(T)}^2 = \vertiii{a}^2 + \|\mathrm{Curl}\,b\|_{L^2(T)}^2.
	\end{align*}
	On the other hand, the constant $m_p$ from \eqref{def:mk} satisfies $$\min_{v_{p+1} \in P_{p+1}(T)} \|q_p - \nabla
	v_{p+1}\|_{L^2(T)}^2 \leq m_p^2\vertiii{\mathrm{curl}\,q_p}_*^2.$$
	Hence, \eqref{ineq:proof-constant-curl-dual-norm} implies
	\begin{align}
		\vertiii{a}^2 \leq (m_p^2C_n^2 - 1)\|\mathrm{Curl}\,b\|_{L^2(T)}^{2}.
		\label{ineq:proof-constant-Grad-a-estimate}
	\end{align}
	The Pythagoras identity $\vertiii{f - a}^2 + \|\mathrm{Curl}\,b\|_{L^2(T)}^2 = \|\nabla f - q_k\|_{L^2(T)}^{2}$, a triangle inequality, the estimate $2(\nabla a, \nabla (f - a))_{L^2(T)} \leq \delta\vertiii{a}^2 + \vertiii{f-a}^2/\delta$, and \eqref{ineq:proof-constant-Grad-a-estimate} show, for all positive parameters $\delta > 0$, that
	\begin{align}
		\vertiii{f - G f}^2 &= \vertiii{f}^2 \,= \vertiii{f - a}^2 + 2(\nabla a, \nabla (f - a))_{L^2(T)} + \vertiii{a}^2\nonumber\\
		&\leq (1 + \delta)\vertiii{f - a}^2 + (1 + 1/\delta)\vertiii{a}^2\nonumber\\
		&\leq \max\{1+\delta,(1+1/\delta)(m_p^2C_n^2-1)\}(\vertiii{f-a}^2 + \|\mathrm{Curl}\,b\|^2)\nonumber\\
		&= \max\{1+\delta, (1+1/\delta)(m_p^2C_n^2-1)\}\|(1 - \Pi_p)\nabla f\|_{L^2(T)}^2.
		\label{ineq:proof-f-Gf-upper-bound}
	\end{align}
	If $m_pC_n > 1$, then $\delta \coloneqq m_p^2 C_n^2 - 1$ leads to $\max\{1+\delta,
	(1+1/\delta)(m_p^2C_n^2-1)\} = m_p^2C_n^2$. If $m_pC_n \leq 1$, then $\inf_{\delta > 0} \max\{1+\delta, (1+1/\delta)(m_p^2C_n^2-1)\} = 1$. This concludes the proof of $C_{\mathrm{st,2}} \leq \max\{1,m_pC_n\}$.
	Notice that $\vertiii{\mathrm{curl}\,q_p}_*^2 = \vertiii{(-\Delta)^{-1} \mathrm{curl}\,q_p}^2 = b(q_p,q_p)$ and
	the orthogonal decomposition $P_p(T;\R^n) = Q_p \oplus \nabla P_{p+1}(T)$ with $\mathrm{curl}\, \nabla P_{p+1}(T)
	\equiv 0$ reveal
	\begin{align}
		m_p^2 = \max_{q_p \in Q_p\setminus\{0\}} \|q_p\|_{L^2(T)}^2/\vertiii{\mathrm{curl}\, q_p}_{*}^2 = \max_{q_p \in Q_p\setminus\{0\}} \mathfrak{a}(q_p,q_p)/\mathfrak{b}(q_p,q_p)
		\label{eq:mk-Rayleigh-quotient}
	\end{align}
	with the bilinear forms $\mathfrak{a}$ and $\mathfrak{b}$ from \eqref{def:a-b}.
	The min-max principle \cite[Sec.\ 8]{BabuskaOsborn1991} and \eqref{eq:mk-Rayleigh-quotient} show that $m_p^2$ is the maximal eigenvalue of \eqref{def:mk-eigenvalue-problem}. This concludes the proof.
\end{proof}
\begin{example}[numerical example]
	\label{ex:C_st}
	\Cref{tab:C_k} displays the computed maximal eigenvalue
	$m_p^2\geq C_{\rm st,2}^2$ of the eigenvalue problem \eqref{def:mk-eigenvalue-problem} for the right-isosceles triangle $T$.
	The right-hand side is approximated by the Courant FEM of polynomial degree $10$ on a uniform triangulation of $T$ with $50721$
	degrees of freedom.
	The lower bounds
	\begin{align*}
		\sup_{f\in P_N(T)}\frac{\vertiii{(1-G)f}}{\|(1-\Pi_p)\nabla f\|_{L^2(T)}}\leq C_{\rm st,2} \quad\text{and}\quad
		\sup_{f\in P_N(T)}\frac{\vertiii{(1-\Pi_{p+1})f}}{\|(1-\Pi_p)\nabla f\|_{L^2(T)}}\leq C_{\rm st,1}
	\end{align*}
	for $C_{\rm st,2}$ and $C_{\rm st,1}$ from \eqref{ineq:stability-estimate-1} are computable Rayleigh quotients and displayed in \Cref{fig:C_p}.
	Computer experiments provide numerical evidence for the convergence of the lower bounds of $C_{\rm st,2}$ to $m_p$
	as $N\to\infty$ and, hence, for $C_{\rm st,2}=m_p$.
	The lower bound of $C_{\rm st, 1}\propto\sqrt{p+1}$ displays
	the expected growth.
\end{example}
Undisplayed numerical experiments suggest that a small minimal interior angle does \emph{not} affect the asymptotic bound of $C_{\rm st,2}$,
but leads to increased growth of $C_{\rm st,1}$ as $p\to\infty$.
We observed $C_{\rm st, 2}=m_p$ and the convergence $C_{\rm st, 2}^2\to 2$ as $p\to \infty$ for different isosceles and
various right triangles, whereas an interior angle $\omega>\pi/2$ has a mild influence on the
	maximal value of $C_{\rm
st, 2}$ as shown for isosceles triangles in \Cref{fig:C_p_angle}.
(Recall that the constants $C_{\rm st, 1}$ and $C_{\rm st, 2}$ are invariant under scaling.)
This leads to our following conjecture in accordance with \Cref{fig:C_p} for any $p\in\mathbb N_0$.
\begin{conjecture}
	For triangles $T$ with maximal interior angle $\omega\leq\pi/2$, $C_{\rm st,2} \leq \sqrt{2}$.
\end{conjecture}

\section{The modified HHO method}\label{sec:HHO}
This section introduces the HHO method and the discrete eigenvalue problem.
\subsection{Triangulation}\label{sec:triangulation}
Let $\Tcal$ be a regular triangulation of $\Omega$ into simplices in the sense of Ciarlet such that $\cup_{T \in \Tcal} T = \overline{\Omega}$.
Given a simplex $T \in \Tcal$ of positive volume $|T| > 0$, let $\Fcal(T)$ denote the set of the $n+1$ hyperfaces of $T$, called sides of $T$. Define the set of all sides $\Fcal = \cup_{T \in \Tcal} \Fcal(T)$ and the set of interior sides $\Fcal(\Omega) = \Fcal\setminus\{F \in \Fcal:F\subset\partial \Omega\}$ in $\Tcal$.
For any interior side $F \in \Fcal(\Omega)$, there exist exactly two simplices $T_+, T_- \in \Tcal$ such that $\partial T_+ \cap \partial T_- = F$. The orientation of the outer normal unit $\nu_F = \nu_{T_+}|_F = -\nu_{T_-}|_F$ along $F$ is fixed and $\omega_F \coloneqq \mathrm{int}(T_+ \cup T_-)$ denotes the side patch of $F$. Let $[v]_F \coloneqq (v|_{T_+})|_F - (v|_{T_-})|_F \in L^1(F)$ denote the jump of $v \in L^1(\omega_F)$ with $v \in W^{1,1}(T_\pm)$ across $F$.
For any boundary side $F \in \Fcal(\partial \Omega) \coloneqq \Fcal\setminus\Fcal(\Omega)$, there exists a unique $T \in \Tcal$ with $F \in \Fcal(T)$.
Then $\omega_F \coloneqq \mathrm{int}(T)$, $\nu_F \coloneqq \nu_T$ is the exterior unit vector of $F \in \Fcal(T)$, and $[v]_F \coloneqq v|_F$. The triangulation $\Tcal$ gives rise to the space $H^1(\Tcal) \coloneqq \{v \in L^2(\Omega): v|_T \in H^1(T)\}$ of piecewise Sobolev functions. The differential operators $\div_\pw$, $\nabla_\pw$, and $\Delta_\pw$ denote the piecewise applications of $\div$, $\nabla$, and $\Delta$ without explicit reference to the triangulation $\Tcal$. 

\subsection{Discrete spaces}

Let $P_p(\Tcal)$, $P_p(\Fcal)$, and $\RT^\pw_p(\Tcal)$ denote the space of piecewise functions with restrictions to $T \in \Tcal$ or $F \in \Fcal$ in $P_p(T)$, $P_p(F)$, and $\RT_p(T)$.
The local mesh-sizes give rise to the piecewise constant function $h_\Tcal \in P_0(\Tcal)$ with $h_\Tcal|_T \equiv h_T = \mathrm{diam}(T)$ in $T \in \Tcal$ and $h_{\max} \coloneqq \|h_\Tcal\|_{L^\infty(\Omega)}$ abbreviates the maximal mesh-size of $\Tcal$.
The $L^2$ projections $\Pi_p : L^1(\Omega) \to P_p(\Tcal)$, $\Pi_\Fcal^p : L^1(\cup \Fcal) \to P_p(\Fcal)$, and $\Pi_{\RT} : L^1(\Omega;\R^n) \to \RT_p^\pw(\Tcal)$ onto $P_p(\Tcal)$, $P_p(\Fcal)$, and $\RT_p^\pw(\Tcal)$ are computed cell-wise.
For vector-valued functions $\tau \in L^1(\Omega;\R^n)$, the $L^2$ projection $\Pi_p$ onto $P_p(\Tcal;\R^n) = P_p(\Tcal)^n$ applies componentwise. The Pythagoras theorem implies the stability of $L^2$ projections, for any $\tau \in L^2(\Omega;\R^n)$ and $v \in L^2(\Omega)$,
\begin{align}
	\|\Pi_\RT \tau\|_{L^2(\Omega)}^2 \leq \|\tau\|_{L^2(\Omega)} \text{ and } \|\Pi_p v\|_{L^2(\Omega)} \leq \|v\|_{L^2(\Omega)}.
	\label{ineq:L2-orthogonality}
\end{align}
The Galerkin projection $G f$ of $f \in H^1(\Tcal)$ is computed cell-wise by \eqref{def:Galerkin} with
\begin{align}
	\|\nabla_\pw(1 - G) f\|_{L^2(\Omega)} = \min_{p_{p+1} \in P_{p+1}(\Tcal)} \|\nabla_\pw(f - p_{p+1})\|_{L^2(\Omega)}.
	\label{ineq:G-best-approximation}
\end{align}
The inclusion $\nabla_\pw P_{p+1}(\Tcal) \subset P_p(\Tcal;\R^n) \subset \RT_p^\pw(\Tcal)$ leads, for any $f \in H^1(\Tcal)$, to
\begin{align}
	\|(1 - \Pi_{\RT}) \nabla_\pw f\|_{L^2(\Omega)} \leq \|(1 - \Pi_p) \nabla_\pw f\|_{L^2(\Omega)} \leq \|\nabla_\pw(1 - G) f\|_{L^2(\Omega)}.
	\label{ineq:inclusion}
\end{align}

\subsection{HHO methodology}\label{sec:HHO-methodology}
Let $V_h \coloneqq P_{p+1}(\Tcal) \times P_p(\Fcal(\Omega))$ denote the ansatz space of the HHO method for $p\in\mathbb
N_0$.
The interior sides $\Fcal(\Omega)$ give rise to the subspace $P_p(\Fcal(\Omega))$ of all $(v_F)_{F \in \Fcal} \in P_p(\Fcal)$ with the convention $v_F \equiv 0$ on any boundary side $F \in \Fcal(\partial \Omega)$ for homogenous boundary conditions. In other words, the notation $v_h \in V_h$ means $v_h = (v_\Tcal,v_\Fcal) = \big((v_T)_{T \in \Tcal},(v_F)_{F \in \Fcal}\big)$ for some $v_\Tcal \in P_{p+1}(\Tcal)$ and $v_\Fcal \in P_p(\Fcal(\Omega))$ with the identification $v_T = v_\Tcal|_T \in P_{p+1}(T)$ and $v_F = v_\Fcal|_F \in P_p(F)$.
Given $v_h = (v_\Tcal,v_\Fcal) \in V_h$, the norm $\|v_h\|_h$ of $v_h$ in $V_h$ from \cite[Eq.~(28)]{DiPietroErnLemaire2014} or \cite[Eq.~(41)]{DiPietroErn2015} reads
\begin{align}\label{def:discrete-norm}
	\|v_h\|_h^2 \coloneqq \|\nabla_\pw v_\Tcal\|^2_{L^2(T)} + \sum_{T \in \Tcal} \sum_{F \in \Fcal(T)} h_F^{-1}\|v_F - v_T\|^2_{L^2(F)}.
\end{align}
The interpolation $\I: V \to V_h$ maps $v\mapsto\I v \coloneqq (\Pi_{p+1} v, \Pi_\Fcal^p v) \in V_h$.\medskip

\noindent\emph{Potential reconstruction.} The potential reconstruction $\PotRec v_h \in P_{p+1}(\Tcal)$ of $v_h = (v_\Tcal,v_\Fcal) \in V_h$ satisfies, for all discrete test functions $\varphi_h \in P_{p+1}(\Tcal)$, that
\begin{align}
	&(\nabla_\pw \PotRec v_h, \nabla_\pw \varphi_h)_{L^2(\Omega)}\nonumber\\
	&\qquad= -(v_\Tcal,\Delta_\pw \varphi_h)_{L^2(\Omega)} + \sum_{F \in \Fcal(\Omega)} \int_F v_F [\nabla_\pw \varphi_h \cdot \nu_F]_F \d{s}.\label{def:R-1}
\end{align}
The bilinear form $(\nabla_\pw \bullet, \nabla_\pw \bullet)_{L^2(\Omega)}$ on the left-hand side of \eqref{def:R-1} defines a scalar product and the right-hand side of \eqref{def:R-1} is a linear functional in the quotient space $P_{p+1}(\Tcal)/P_0(\Tcal)$.
The Riesz representation $\PotRec v_h \in P_{p+1}(\Tcal)$ of this linear functional in $P_{p+1}(\Tcal)/P_0(\Tcal)$ is selected by
\begin{align}
	\Pi_0 \PotRec v_h = \Pi_0 v_\Tcal.
	\label{def:R-2}
\end{align}
The unique solution $\PotRec v_h \in P_{p+1}(\Tcal)$ to \eqref{def:R-1}--\eqref{def:R-2} defines the potential reconstruction operator $\PotRec : V_h \to P_{p+1}(\Tcal)$.\medskip

\noindent\emph{Gradient reconstruction.} The gradient reconstruction $\GrRec v_h \in \RT_p^\pw(\Tcal)$ of $v_h =
(v_\Tcal,v_\Fcal) \in V_h$ satisfies, for all discrete test functions $\phi_h \in \RT_p^\pw(\Tcal)$, that
\begin{align}
	(\GrRec v_h, \phi_h)_{L^2(\Omega)} = -(v_\Tcal,\div_\pw \phi_h)_{L^2(\Omega)} + \sum_{F \in \Fcal(\Omega)} \int_F v_F [\phi_h \cdot \nu_F]_F \d{s}.
	\label{def:gradient-reconstruction}
\end{align}
In other words, $\GrRec v_h$ is the Riesz representation of the linear functional on the right-hand side of \eqref{def:gradient-reconstruction} in the Hilbert space $\RT_p^\pw(\Tcal)$ endowed with the $L^2$ scalar product.
Since $\nabla_\pw P_{p+1}(\Tcal) \subset \RT_p^\pw(\Tcal)$, \eqref{def:gradient-reconstruction} implies the $L^2$ orthogonality $\GrRec v_h - \PotRec v_h \perp \nabla_\pw P_{p+1}(\Tcal)$.
The following lemma recalls the commutativity of $\GrRec$ and $\PotRec$
\cite{DiPietroErnLemaire2014,DiPietroErn2015,AbbasErnPignet2018,DiPietroTittarelli2018}.
The Galerkin projection $G$ is
defined in \eqref{def:Galerkin}.
\begin{lemma}[commutativity]
	\label{lem:properties-reconstruction-operators} Any $v \in V$ satisfies $\GrRec \I v = \Pi_{\RT} \nabla v$ and
$\PotRec \I v = G v$.\qed
\end{lemma} 

\subsection{Discrete eigenvalue problem}\label{sec:discrete-problem}
Given positive constants $0 < \alpha < 1$ and $0 < \beta < \infty$,  recall $a_h$ and $b_h$ from
\eqref{def:ah}--\eqref{def:bh}. Notice, for any $u_h, v_h \in V_h$, that
\begin{align}
	a_h(u_h, v_h) = ((1 - \alpha) \GrRec u_h + \alpha \Pi_p \GrRec u_h, \GrRec v_h)_{L^2(\Omega)} + \beta(h_\Tcal^{-2}S u_h, S v_h)_{L^2(\Omega)}.
	\label{eq:ah-alternative-form}
\end{align}
The discrete problem seeks a discrete eigenpair $(\lambda_h, u_h) \in \R_+ \times V_h\setminus\{0\}$ such that
\begin{align}
	a_h(u_h,v_h) = \lambda_h b_h(u_h, v_h) \text{ for all } v_h \in V_h.
	\label{pr:discrete-problem}
\end{align}
\begin{lemma}[discrete norm]\label{lem:discrete-norm}
	The bilinear form $a_h: V_h \times V_h \to \R$ is a scalar product in $V_h$. The induced norm $\|\bullet\|_{a,h} \coloneqq a_h(\bullet,\bullet)^{1/2} \approx \|\bullet\|_h$ is equivalent to the discrete norm $\|\bullet\|_h$ from \eqref{def:discrete-norm}.
\end{lemma}
\begin{proof}
	The equivalence $\|\bullet\|_{a,h} \approx \|v_h\|_h$ for all $v_h \in V_h$ is proven in \cite[Lemma
	3.5]{CarstensenErnPuttkammer2021} for the stabilization $\beta(\nabla_\pw S u_h, \nabla_\pw S v_h)_{L^2(\Omega)}$
	instead of $\beta(h_\Tcal^{-2} S u_h, S v_h)_{L^2(\Omega)}$ in the definition \eqref{def:ah} of $a_h$. Since the two
	stabilizations are locally equivalent, this leads to the assertion.
\end{proof}
The discrete eigenvalue problem \eqref{pr:discrete-problem} gives rise to the symmetric generalized algebraic eigenvalue problem
\begin{align}
	\begin{pmatrix}
		A_{\Tcal \Tcal} & A_{\Tcal \Fcal}\\ A_{\Fcal \Tcal} & A_{\Fcal \Fcal}
	\end{pmatrix}
	\begin{pmatrix}
		x_\Tcal\\ x_\Fcal
	\end{pmatrix} = \lambda_h
	\begin{pmatrix}
		B_{\Tcal\Tcal} & 0\\ 0 & 0
	\end{pmatrix}
	\begin{pmatrix}
		x_\Tcal\\ x_\Fcal
	\end{pmatrix}.
	\label{eq:algebraic-eigenvalue-problem}
\end{align}
The application of the Schur complement as in \cite[Section 3.3]{CarstensenErnPuttkammer2021} leads to the algebraic
eigenvalue problem $(A_{\Tcal\Tcal} - A_{\Tcal\Fcal}A_{\Fcal\Fcal}^{-1}A_{\Fcal\Tcal}) x_\Tcal = \lambda_h
B_{\Tcal\Tcal} x_\Tcal$. Hence, \eqref{eq:algebraic-eigenvalue-problem} provides $N \coloneqq \dim
P_{p+1}(\Tcal)=|\Tcal|\bigl(\begin{smallmatrix}
	p+1+n\\n
\end{smallmatrix}\bigr)$ positive discrete eigenvalues $0 < \lambda_h(1) \leq \lambda_h(2) \leq \dots \leq \lambda_h(N) < \infty$; all other eigenvalues $\lambda_h(j) \coloneqq \infty$ for $j > N$ are infinity.

\section{Lower eigenvalue bounds}\label{sec:GLB}
This section establishes the sufficient conditions on the parameters $\alpha, \beta$ in \eqref{eq:cond2} such that the HHO method from \eqref{pr:discrete-problem} provides direct GLB.
Let $\lambda$ (resp.~$\lambda_h$) denote the $j$-th continuous (resp.~discrete) eigenvalue of \eqref{pr:eigenvalue-problem} (resp. \eqref{pr:discrete-problem}) for 
fixed $j \in \mathbb{N}$.
Recall $0 < \alpha < 1$, $0 < \beta < \infty$, and the constant $\sigma_2$ from \eqref{ineq:stability-2}.

\begin{theorem}[GLB]\label{thm:leb}
	If $\sigma_2^2\max\{\beta, h_{\max}^2\min\{\lambda_h,\lambda\}\} \leq \alpha$, then $\lambda_h \leq \lambda$.
\end{theorem}
\begin{remark}[GLB for $j > N$]
	The number $j\in\mathbb N$ in the theorem can be larger than the dimension $N$.
	Then $\alpha<\sigma_2^2\lambda h_{\max}^2$ follows.
	In other words $\lambda(N+1)>\alpha\sigma_2^{-2}h_{\max}^{-2}$ is an a priori bound for the exact eigenvalue $\lambda(N+1)$ for free.
\end{remark}

\begin{proof}[Proof of \Cref{thm:leb}]
	The proof applies the key arguments  from  \cite[Theorem 4.1]{CarstensenErnPuttkammer2021}, but then reflects a different stabilization. This enables a different 
	sufficient condition in the theorem
	with a more appropriate precise arrangement of
	the constants. (In fact, $G$ in \eqref{ineq:stability-1}-\eqref{ineq:stability-2} is replaced by $\Pi_{p+1}$ in  \cite{CarstensenErnPuttkammer2021}, whence $C_{\mathrm{st,2}}$ in this 
	paper is not larger than $C_{\mathrm{st},1}$ in \cite{CarstensenErnPuttkammer2021} and $\kappa$ from   \cite{CarstensenErnPuttkammer2021} is bounded by $\sigma_2$ from 
	\eqref{ineq:stability-2}.) Besides those differences, the first steps in the proof  are very analogous and adopted for brevity.
	
	Observe carefully that, in the beginning,   $\sigma_2^2 h_{\max}^2\min\{\lambda_h,\lambda\} \leq \alpha$ does not
	immediately imply that $0<\lambda_h\leq\infty$ is finite. 

	\medskip
	
	\noindent\emph{Step 1: Reduction to $h_{\max}^2 \sigma_2^2 \lambda < 1$.}  If $h_{\max}^2 \sigma_2^2 \lambda \geq 1$, then $h_{\max}^2 \sigma_2^2 \lambda_h \leq \alpha < 1 \leq h_{\max}^2 \sigma_2^2 \lambda$, 	whence $\lambda_h$ is finite and $\lambda_h \leq \lambda$. The remaining parts of this proof therefore 
	assume $h_{\max}^2 \sigma_2^2 \lambda < 1$.
	
	\medskip
	
	\noindent\emph{Step 2: The first $j$ exact and pairwise orthonormal 
		eigenfunctions $\phi_1, \dots, \phi_j \in V$ of \eqref{pr:eigenvalue-problem} satisfy that $\Pi_{p+1} \phi_1, \dots, \Pi_{p+1} \phi_j \in P_{p+1}(\Tcal)$ are linear independent.}
	The proof follows the lines of Step 2 in the proof of \cite[Theorem 4.1]{CarstensenErnPuttkammer2021} (with %
	$\delta \coloneqq \sigma_2 h_{\max}$).
	
	\medskip
	
	\noindent\emph{Step 3: There exists $\phi \in \mathrm{span}\{\phi_1,\dots,\phi_j\}$ with $\|\phi\|_{L^2(\Omega)} = 1$, $\|\nabla \phi\|^2_{L^2(\Omega)} \leq \lambda$, and}
	\begin{align}
		0<\lambda_h (1 - \|(1 - \Pi_{p+1}) \phi\|_{L^2(\Omega)}^2) \leq a_h(\I \phi, \I \phi).
		\label{ineq:leb-step-3}
	\end{align}
	The proof follows the lines of Step~3 in the proof of \cite[Theorem 4.1]{CarstensenErnPuttkammer2021} and considers the $\min$-$\max$ principle for the 
	algebraic eigenvalue problem \eqref{eq:algebraic-eigenvalue-problem} with the $j$-dimensional subspace spanned by  $\I\phi_1, \dots, \I\phi_j \in V_h$.
	It is the linear independence of  $\Pi_{p+1} \phi_1, \dots, \Pi_{p+1} \phi_j \in P_{p+1}(\Tcal)$ that guarantees 
	$j \leq N = \dim P_{p+1}(\Tcal)$ and 
	that the algebraic eigenvalue problem 
	\eqref{eq:algebraic-eigenvalue-problem} has at least $j$ finite eigenvalues; whence $\lambda_h=\lambda_h(j)<\infty$. The bound of $\lambda_h$ in the 
	$\min$-$\max$ principle by some maximizer $v_h\coloneqq\I\phi$  of the Rayleigh quotient  in  $ \mathrm{span}\{\I\phi_1, \dots, \I\phi_j \} \subset  V_h$ is rewritten as
	\[
		\lambda_h b_h(\I \phi, \I \phi) \leq a_h(\I \phi, \I \phi)<\infty 
	\] 
	for
	$\phi \in \mathrm{span}\{\phi_1,\dots,\phi_j\}$ with $\|\phi\|_{L^2(\Omega)} = 1$ and $\|\nabla \phi\|^2_{L^2(\Omega)} \leq \lambda$.
	It  follows from Step~2 that $ b_h(\I \phi, \I \phi) =\| \Pi_{p+1}\phi\|_{L^2(\Omega)}^2  >0$ cannot vanish. 
	This and the Pythagoras theorem $\|\Pi_{p+1} \phi\|_{L^2(\Omega)}^2 =1 - \|(1 - \Pi_{p+1})\phi\|_{L^2(\Omega)}^2>0 $ 
	(recall $1=\|\phi\|^2_{L^2(\Omega)} $)
	conclude the proof of \eqref{ineq:leb-step-3}.
	
	\medskip
	
	\noindent\emph{Step 4: First lower bound for $\lambda - \lambda_h$ under the assumption $\beta \sigma_2^2 \le \alpha $.}
	The commutativity $\GrRec \I \phi = \Pi_{\RT} \nabla \phi$ from \Cref{lem:properties-reconstruction-operators}.a 
	and $(1 - \alpha) \GrRec u_h + \alpha \Pi_p \GrRec u_h = (1 - \alpha)(1 - \Pi_p)\GrRec u_h + \Pi_p \GrRec u_h$ for
	$u_h = \I \phi$ 
	prove that $a_h(\I \phi, \I \phi)$ in  \eqref{eq:ah-alternative-form} is equal to
	\begin{align}
		(1 - \alpha)\|(1 - \Pi_p)\Pi_\RT \nabla \phi\|_{L^2(\Omega)}^2 + \|\Pi_p \Pi_\RT \nabla \phi\|_{L^2(\Omega)}^2 + \beta\|h_\Tcal^{-1} S \I \phi\|_{L^2(\Omega)}^2.
		\label{ineq:proof-LEB-step-4}
	\end{align}
The identity $\|(1 - \Pi_p) \Pi_{\RT} \nabla \phi\|_{L^2(\Omega)} = \|\Pi_\RT(1 - \Pi_p) \nabla \phi\|_{L^2(\Omega)}$
follows from the inclusion
	$P_p(\Tcal;\R^n) \subset \RT_p^\pw(\Tcal)$ and 
$\Pi_p \Pi_{\RT} \nabla \phi = \Pi_p \nabla \phi = \Pi_\RT \Pi_p \phi$.
	This, \eqref{ineq:proof-LEB-step-4}, and $\|\Pi_\RT(1 - \Pi_p) \nabla \phi\|_{L^2(\Omega)} \leq \|(1 - \Pi_p) \nabla \phi\|_{L^2(\Omega)}$ from \eqref{ineq:L2-orthogonality} lead to
	\begin{align}
		a_h(\I \phi, \I \phi) \leq \|\Pi_p \nabla \phi\|_{L^2(\Omega)}^2 + (1 - \alpha)\|(1 - \Pi_p) \nabla \phi\|_{L^2(\Omega)}^2 + \beta\|h_\Tcal^{-1} S \I \phi\|_{L^2(\Omega)}^2.
		\label{ineq:proof-LEB-bound-a_h(Iphi,Iphi)}
	\end{align}
	The Pythagoras theorem and $\|\nabla \phi\|_{L^2(\Omega)}^2 \leq \lambda$ prove
	\begin{align}
		\|\Pi_p \nabla \phi\|_{L^2(\Omega)}^2 \leq \lambda - \|(1 - \Pi_p) \nabla \phi\|^2_{L^2(\Omega)}.
	\end{align}
	Recall $S \I \phi = \Pi_{p+1} \phi - \PotRec \I \phi = \Pi_{p+1} \phi - G \phi$ from \Cref{lem:properties-reconstruction-operators}.b. 
	The piecewise mesh-size function $h_\Tcal$  does \emph{not} interfere with the projection $\Pi_{p+1}$ and so the Pythagoras theorem reads
	\begin{align}
		\|h_\Tcal^{-1} S \I \phi\|_{L^2(\Omega)}^2 = \|h_\Tcal^{-1}(1 - G) \phi\|_{L^2(\Omega)}^2 - \|h_\Tcal^{-1}(1 - \Pi_{p+1}) \phi\|_{L^2(\Omega)}^2.
		\label{ineq:leb-step-4-1}
	\end{align}
	The combination of \eqref{ineq:leb-step-3} with \eqref{ineq:proof-LEB-bound-a_h(Iphi,Iphi)}--\eqref{ineq:leb-step-4-1} results in 
	\begin{align*}
		&-\lambda_h \|(1 - \Pi_{p+1}) \phi\|_{L^2(\Omega)}^2 + \alpha \|(1 - \Pi_p)\nabla \phi\|_{L^2(\Omega)}^2\nonumber\\
		&\qquad + \beta \|h_\Tcal^{-1}(1 - \Pi_{p+1}) \phi\|_{L^2(\Omega)}^2 - \beta\|h_\Tcal^{-1}(1 - G) \phi\|_{L^2(\Omega)}^2 \leq \lambda - \lambda_h.
	\end{align*}
	This, the stability estimate \eqref{ineq:stability-2}, and $h_{\max}^{-1} \leq h_\Tcal^{-1}$ in $\Omega$ imply
	\begin{align}
		&(\beta/h_{\max}^2 - \lambda_h)\|(1 - \Pi_{p+1}) \phi\|_{L^2(\Omega)}^2\nonumber\\
		&\qquad\qquad + (\alpha - \beta \sigma_2^2)\|(1 - \Pi_p) \nabla \phi\|_{L^2(\Omega)}^2 \leq \lambda - \lambda_h.
		\label{ineq:leb-step-4}
	\end{align}
	Recall $\|(1 - \Pi_{p+1}) \phi\|_{L^2(\Omega)}^2 \,\leq \|(1 - G) \phi\|_{L^2(\Omega)}^2 \leq \sigma_2^2 h_{\max}^2\|(1 - \Pi_p)\nabla \phi\|_{L^2(\Omega)}^2$ from the best approximation property of $\Pi_{p+1}$ and \eqref{ineq:stability-2}
	as well as $\alpha - \beta \sigma_2^2 \geq 0$ from the  assumptions. Consequently, 
	the left-hand side of \eqref{ineq:leb-step-4} is greater than or equal to $\|(1 - \Pi_{p+1})\phi\|_{L^2(\Omega)}^2$ times 
	\begin{align*}
		(\beta/h_{\max}^2 - \lambda_h + (\alpha - \beta \sigma_2^2)/(\sigma_2^2 h_{\max}^2))=\alpha\sigma_2^{-2} h_{\max}^{-2} - \lambda_h.
	\end{align*}
	In conclusion, $0\le \|(1 - \Pi_{p+1}) \phi\|_{L^2(\Omega)}<1 $ (from the end of Step~3) and 
	\begin{align}
		(\alpha\sigma_2^{-2} h_{\max}^{-2} - \lambda_h)\|(1 - \Pi_{p+1}) \phi\|_{L^2(\Omega)}^2 \leq \lambda - \lambda_h.
		\label{ineq:leb-step-4a}
	\end{align}
	
	\medskip
	
	\emph{Step 5: Finish of the proof.}
	After the reduction to $h_{\max}^2 \sigma_2^2 \lambda < 1$, the 
	above Steps~2--4 of the proof have utilized  $\beta \sigma_2^2 \le \alpha $, but they 
	carefully avoided any assumption on $\lambda$ and $\lambda_h$, although it is supposed that 
	$\sigma_2^2 h_{\max}^2\min\{\lambda_h,\lambda\} \leq \alpha$. In case that 
	$\sigma_2^2 h_{\max}^2\lambda_h \leq \alpha$, the assertion $0\le  \lambda - \lambda_h$ follows immediately from \eqref{ineq:leb-step-4a}.
	In the remaining case  $\sigma_2^2 h_{\max}^2\lambda \leq \alpha$, the pre-factor in the left-hand side of \eqref{ineq:leb-step-4a}
	has the lower bound $\lambda - \lambda_h \le \alpha\sigma_2^{-2} h_{\max}^{-2} - \lambda_h$. Therefore \eqref{ineq:leb-step-4a} implies 
	\[
		(\lambda - \lambda_h) \|(1 - \Pi_{p+1}) \phi\|_{L^2(\Omega)}^2\le \lambda - \lambda_h.
	\]
	Recall that 
	$0\le \|(1 - \Pi_{p+1}) \phi\|_{L^2(\Omega)}<1$ from Step 4 to see that the last displayed estimate  is impossible unless $0 \leq \lambda - \lambda_h$.
	\qedhere
\end{proof}

\section{A priori error analysis}\label{sec:a-priori}
The Babu\v{s}ka-Osborn theory \cite{BabuskaOsborn1991} for the spectral approximation of compact self\-adjoint operators leads to a~priori convergence rates for the approximation of $\lambda$ and of $u$ in the energy norm \cite{CaloCicuttinDengErn2019,CarstensenErnPuttkammer2021}.
This section establishes the quasi-best approximation estimate \eqref{ineq:quasi-optimal-eigenvalue} for a simple
eigenvalue $\lambda$, that eventually allows for the absorption of higher-order terms in the a~posteriori error analysis
of \Cref{sec:a-posteriori}.

Throughout the remaining parts of this paper, suppose that $\beta \leq \alpha/\sigma_2^2$ with $\sigma_2$ from \eqref{ineq:stability-2}.
Let $\lambda \coloneqq \lambda(j)$ be a simple eigenvalue of \eqref{pr:eigenvalue-problem} with
the corresponding eigenfunction $u \coloneqq u(j) \in V$. Let $(\lambda_h,u_h) \coloneqq (\lambda_h(j), u_h(j))$ denote the $j$-th discrete
eigenpair of \eqref{pr:discrete-problem} with $u_h = (u_\Tcal,u_\Fcal) \in V_h$, $\|u\|_{L^2(\Omega)} = 1 =
\|u_\Tcal\|_{L^2(\Omega)}$, and $(u,u_\Tcal)_{L^2(\Omega)} \geq 0$. Recall that $0 < s \leq 1$ denotes the minimum of the index of elliptic regularity and one.

\begin{theorem}[a~priori]\label{thm:L2-error-control}
	If $h_{\max}$ is sufficiently small, then \eqref{ineq:quasi-optimal-eigenvalue} holds.
	The constant $\cnst{cnst:best-approximation}$ exclusively depends on $p$, $n$, $\Omega$, and the shape regularity of $\mathcal{T}$.
\end{theorem}
The following lemmas precede the proof of \Cref{thm:L2-error-control}. The first one recalls the enriching operator from
\cite{ErnZanotti2020} and adds the estimate \eqref{ineq:J-best-approximation}.
Recall the induced discrete norm $\|\bullet\|_{a,h} \coloneqq a_h(\bullet,\bullet)^{1/2}$ from \Cref{lem:discrete-norm}.
\begin{lemma}[enriching operator]\label{lem:enriching-operator}
	There exists a linear bounded operator $\mathrm{J} : V_h \to V$ that is a right-inverse of $\mathrm{I}$, i.e., $v_h=\I
	\mathrm{J} v_h=(\Pi_{p+1}\mathrm Jv_h, \Pi_\Fcal^p\mathrm Jv_h)$ for all $v_h \in V_h$, and stable in the sense that $\|\nabla \mathrm{J} v_h\|_{L^2(\Omega)}
	\leq \|\mathrm{J}\|\|v_h\|_{a,h}$ with $\|\mathrm{J}\|\lesssim 1$. Any $v \in V$ satisfies
	\begin{align}
		\|\nabla(v - \mathrm{J} \I v)\|_{L^2(\Omega)} \leq \cnst{cnst:enrichment-operator} \min_{v_{p+1} \in P_{p+1}(\Tcal)} \|\nabla_\pw(v - v_{p+1})\|_{L^2(\Omega)}.
		\label{ineq:J-best-approximation}
	\end{align}
	The constants $\|\mathrm{J}\|$ and $\newcnst\label{cnst:enrichment-operator}$ solely depend on $p$, $n$, and the shape regularity of $\mathcal{T}$.
\end{lemma}
\begin{proof}
	The construction of the enriching operator $\mathrm{J}:V_h\to V$ in spirit of \cite{VeeserZanottiII} involves standard averaging and
	bubble-function techniques from \cite{Verfuerth1986} and is explained in \cite[Section 4.3]{ErnZanotti2020} for
	a related HHO method without the proof of \eqref{ineq:J-best-approximation}.
	Notice that $\mathrm{J}$ from \cite{ErnZanotti2020} (called stabilized bubble smoother
	$E_\mathrm{H}$ therein) only satisfies $\mathrm{J} v_h - v_\Tcal \perp P_{p-1}(\Tcal)$
	for any given $v_h=(v_\Tcal, v_\Fcal)\in V_h$.
	However, a straight-forward modification of \cite[Eq.~(4.16)]{ErnZanotti2020}
	(in the notation of \cite{ErnZanotti2020}, $\mathcal{B}_K v_\mathcal{M} \in P_{p+1}(K)$ should be defined by
	equation (4.16)
	therein for all $q \in P_{p+1}(K)$) immediately provides a right-inverse $\mathrm{J}$ of $\mathrm{I}$.
	The arguments from \cite[Propositions 4.5 and 4.7]{ErnZanotti2020} apply and lead to the stability of $\mathrm{J}$ with
	respect to the equivalent discrete norm $\|\bullet\|_h\approx\|\bullet\|_{a,h}$ from \Cref{lem:discrete-norm}.

	It remains to prove \eqref{ineq:J-best-approximation} which is well-known for the Crouzeix-Raviart finite element
	method with an appropriate interpolation $\mathrm{I}$ and the conforming companion $\mathrm{J}$ from
	\cite[Proposition 2.3]{CarstensenGallistlSchedensack2015} for $n=2$ and from \cite[Section 5.8]{carstensen_how_2020} for $n=3$.
	Given any $v_h\in V_h$, let $\mathcal{A} v_h \in S^{p+1}_0(\Tcal)\coloneqq P_{p+1}(\Tcal)\cap H^1_0(\Omega)$ denote the nodal average of $\PotRec v_h$, cf.~\cite[Eq.~(4.24)]{ErnZanotti2020}.
	With \cite[Eq.~(4.18)]{ErnZanotti2020} and with the above modification in \cite[Eq.~(4.16)]{ErnZanotti2020}, the bubble smoother $\mathcal B:L^2(\Omega)\times
	L^2(\bigcup\Fcal)\to H^1_0(\Omega)$ from
	\cite[Proposition 4.6]{ErnZanotti2020} satisfies, for $(v_{\mathcal{M}},v_\Sigma)\in L^2(\Omega)\times
	L^2(\bigcup\Fcal)$, the stability estimate
	\begin{align}\label{ineq:J-stability-average-bubble}
		\|\nabla\mathcal{B}(v_{\mathcal{M}}, v_{\Sigma})\|_{L^2(\Omega)}^2\lesssim
		\|h_\Tcal^{-1}\Pi_{p+1}v_{\mathcal{M}}\|_{L^2(\Omega)}^2 + \sum_{T\in\Tcal}\sum^{}_{F\in\Fcal(T)}
		h_F^{-1}\|\Pi_F^p v_\Sigma\|_{L^2(F)}^2
	\end{align}
	with the $L^2$ projection $\Pi_F^p$ onto $P_p(F)$ for all faces $F \in \Fcal$.
	A triangle inequality, the stability of $\Pi_F^p$ on a face $F$, and a discrete
	trace inequality show $\|\Pi_F^p (v_F - \mathcal{A} v_h)\|_{L^2(F)} \leq \|\Pi_F^p (v_F - (\PotRec v_h)|_T)\|_{L^2(F)} +
	h_F^{-1/2}\|\PotRec v_h - \mathcal{A} v_h\|_{L^2(T)}$ for all $F \in \Fcal(T)$ and $T \in \Tcal$. This, a triangle
	inequality for $\mathrm J\coloneqq \mathcal{A}+\mathcal{B}(1-\mathcal{A})$, \eqref{ineq:J-stability-average-bubble}, and the second inequality on \cite[p.~2180]{ErnZanotti2020} result in
	\begin{align}
		&\|\nabla_\pw(\PotRec v_h - \mathrm{J} v_h)\|_{L^2(\Omega)}^2 \lesssim \sum_{F \in \Fcal} h_F^{-1}\|[\PotRec v_h]_F\|_{L^2(\Omega)}^2\nonumber\\
		&\qquad + \|h_\Tcal^{-1} (v_\Tcal - \PotRec v_h)\|_{L^2(\Omega)}^2 + \sum_{T \in \Tcal} \sum_{F \in \Fcal(T)}
		h_F^{-1}\|\Pi_F^p (v_F - (\PotRec v_h)|_T)\|_{L^2(F)}^2.
		\label{ineq:J-stability-jump-stabilization}
	\end{align}
	Given $v \in V$,
	the stability of the $L^2$ projections $\Pi_{p+1}$ and $\Pi_F^p$ from \eqref{ineq:L2-orthogonality} prove
	$\|\Pi_{p+1}(v - \PotRec \I v)\|_{L^2(T)} \leq \|v - \PotRec \I v\|_{L^2(T)}$ and $\|\Pi_F^p (v - \PotRec (\I v)|_T)\|_{L^2(F)} \leq \|v - (\PotRec \I v)|_T\|_{L^2(F)}$ for all $T \in \Tcal$ and $F \in \Fcal(T)$.
	Given an interior side $F = T_+ \cap T_- \in \Fcal(\Omega)$ for $T_\pm \in \Tcal$, the triangle inequality shows 
	\begin{align*}
		\|[\PotRec \I v]_F\|_{L^2(F)} = \|[\PotRec \I v - v]_F\|_{L^2(F)} \leq \|(v - \PotRec \I v)|_{T_+}\|_{L^2(F)} + \|(v - \PotRec \I v)|_{T_-}\|_{L^2(F)}.
	\end{align*}
	For boundary sides $F \in \Fcal(\partial \Omega)$, it holds $\|[\PotRec \I v]_F\|_{L^2(F)} = \|v - \PotRec \I v\|_{L^2(F)}$.
	The choice $v_h \coloneqq \I v$ in \eqref{ineq:J-stability-jump-stabilization}, the aforementioned inequalities, the trace inequality, and the piecewise application of the Poincar\'e inequality imply $\|\nabla_\pw(\PotRec \I v - \mathrm{J} \I v)\|_{L^2(\Omega)} \lesssim \|\nabla_\pw(v - \PotRec \I v)\|_{L^2(\Omega)}$.
	This, the triangle inequality 
	\begin{align*}
		\|\nabla(v - \mathrm{J} \I v)\|_{L^2(\Omega)} \leq \|\nabla_\pw(v - \PotRec \I v)\|_{L^2(\Omega)} +
		\|\nabla_\pw(\PotRec \I v - \mathrm{J} \I v)\|_{L^2(\Omega)},
	\end{align*}
	and the $L^2$ orthogonality $\nabla_\pw(v - \PotRec \I v) \perp \nabla_\pw P_{p+1}(\Tcal)$ conclude the proof of
	\eqref{ineq:J-best-approximation}.
\end{proof}
\noindent The second lemma proves quasi-best approximation estimates for a source problem.
\begin{lemma}[best-approximation]\label{lem:best-approximation}
	Given $f \in L^2(\Omega)$, let $\widetilde{u} \in V$ solve $-\Delta \widetilde{u} = f$ in $\Omega$.
	The solution $\widetilde{u}_h = (\widetilde{u}_\Tcal, \widetilde{u}_\Fcal) \in V_h$ to
	\begin{align}
		a_h(\widetilde{u}_h, v_h) = (f, v_\Tcal)_{L^2(\Omega)} \quad\text{for all } v_h = (v_\Tcal,v_\Fcal) \in V_h
		\label{def:auxiliary-problem}
	\end{align}
	and the data oscillation $\mathrm{osc}(f,\Tcal) \coloneqq \|h_\Tcal(1 - \Pi_{p+1})f\|_{L^2(\Omega)}$ satisfy
	\begin{align}
		\cnst{cnst:best-approximation-source}^{-1}\|\I \widetilde{u} - \widetilde{u}_h\|_{a,h} \leq \min_{v_{p+1} \in P_{p+1}(\Tcal)} \|\nabla_\pw(\widetilde{u} - v_{p+1})\|_{L^2(\Omega)} + \mathrm{osc}(f,\Tcal)
		\label{ineq:best-approximation}
	\end{align}
	with the constant $\newcnst\label{cnst:best-approximation-source}\coloneqq\max\{\|\mathrm{J}\| + (\alpha^2/(1-\alpha)+ \beta C_P^2)^{1/2}, \|\mathrm{J}\|C_P\}$.
\end{lemma}
\begin{proof}
	Throughout this proof, abbreviate $\widetilde{e}_h \coloneqq \I \widetilde{u} - \widetilde{u}_h \in V_h$. 
	Since $\Pi_{p+1}u - \widetilde{u}_\Tcal=\Pi_{p+1}\mathrm J\widetilde{e}_h\in
	P_{p+1}(\Tcal)$ by \Cref{lem:enriching-operator}, the discrete problem \eqref{def:auxiliary-problem} shows
	\begin{align}\label{eqn:a_h_5.3}a_h(\widetilde{u}_h, \widetilde{e}_h) =(f, \Pi_{p+1}\mathrm J\widetilde{e}_h)_{L^2(\Omega)}.\end{align}
	The commutativity $\Pi_{\RT} \nabla v = \GrRec \I v$ for $v\in V$ from
	\Cref{lem:properties-reconstruction-operators} enters this proof in two ways.
	First, it verifies $\Pi_p \GrRec \I \widetilde{u} = \Pi_p \nabla
	\widetilde{u}$ with $v\coloneqq \widetilde u$ so that \eqref{eq:ah-alternative-form} reads
	\begin{align}
		a_h(\I\widetilde u, \widetilde{e}_h) &= (\nabla \widetilde{u}-\alpha(1 - \Pi_p)\nabla \widetilde{u}, \GrRec
		\widetilde{e}_h)_{L^2(\Omega)} + \beta(h_\Tcal^{-2} S \I\widetilde{u}, S \widetilde{e}_h)_{L^2(\Omega)}. \label{eq:proof-best-approximation-split}
	\end{align}
Second, for $v\coloneqq \mathrm J\widetilde e_h$, the resulting
	$L^2$ orthogonality $\nabla \mathrm{J} \widetilde{e}_h - \GrRec \widetilde{e}_h \perp \RT_p^\pw(\Tcal)$ to the
	piecewise Raviart-Thomas functions $\RT_p^\pw(\Tcal)$ provides
	\begin{align}
		(\nabla \widetilde{u}, \GrRec \widetilde{e}_h)_{L^2(\Omega)}
		= -((1 - \Pi_\RT) \nabla \widetilde{u},\nabla \mathrm{J} \widetilde{e}_h)_{L^2(\Omega)} + (\nabla \widetilde{u}, \nabla \mathrm{J} \widetilde{e}_h)_{L^2(\Omega)}.\nonumber
	\end{align}
	Since $\widetilde u\in V$ solves $-\Delta \widetilde{u} = f$ in $\Omega$, this and
	\eqref{eqn:a_h_5.3}--\eqref{eq:proof-best-approximation-split} verify
	\begin{align}
		\|\widetilde{e}_h\|^2_{a,h} &= (f, (1 - \Pi_{p+1}) \mathrm{J} \widetilde{e}_h)_{L^2(\Omega)} - ((1-\Pi_{\RT})\nabla \widetilde{u}, \nabla \mathrm{J} \widetilde{e}_h)_{L^2(\Omega)}\nonumber\\
		&\qquad - \alpha((1 - \Pi_p)\nabla \widetilde{u},\GrRec \widetilde{e}_h)_{L^2(\Omega)} + \beta(h_\Tcal^{-2} S \mathrm{I} \widetilde{u}, S \widetilde{e}_h)_{L^2(\Omega)}.
		\label{ineq:proof-best-approximation}
	\end{align}
	The choice $\phi \coloneqq \widetilde{u}$ in \eqref{ineq:leb-step-4-1} implies
	$\|h_\Tcal^{-1} S\I \widetilde{u}\|_{L^2(\Omega)} \leq \|h_\Tcal^{-1}(\widetilde{u} - G \widetilde{u})\|_{L^2(\Omega)}$ with the Galerkin projection $G$ from \eqref{def:Galerkin}. Hence, the Poincar\'e inequality
	shows
	\begin{align}
		\|h_\Tcal^{-1} S \mathrm{I} \widetilde{u}\|_{L^2(\Omega)} \leq C_P \|\nabla_\pw(\widetilde{u} - G \widetilde{u})\|_{L^2(\Omega)}.
		\label{ineq:stabilization-best-approximation}
	\end{align}
	A Cauchy and a piecewise application of the Poincar\'e inequality reveal
	\begin{align}
		(f, (1 - \Pi_{p+1}) \mathrm{J} \widetilde{e}_h)_{L^2(\Omega)} \leq C_P \mathrm{osc}(f,\Tcal)\|\nabla \mathrm{J} \widetilde{e}_h\|_{L^2(\Omega)}.
		\label{ineq:h.o.t}
	\end{align}
	The combination of \eqref{ineq:proof-best-approximation}--\eqref{ineq:h.o.t} with a Cauchy inequality provides
	\begin{align*}
		&\|\widetilde{e}_h\|_{a,h}^2 \leq \big(C_P\mathrm{osc}(f,\Tcal) + \|(1 - \Pi_{\RT}) \nabla \widetilde{u}\|_{L^2(\Omega)}\big)\|\nabla \mathrm{J} \widetilde{e}_h\|_{L^2(\Omega)}\\
		&\quad+ \alpha\|(1 - \Pi_p) \nabla \widetilde{u}\|_{L^2(\Omega)}\|\GrRec \widetilde{e}_h\|_{L^2(\Omega)} + \beta
		C_P\|\nabla_\pw(\widetilde{u} - G \widetilde{u})\|_{L^2(\Omega)}\|h_\Tcal^{-1}S \widetilde{e}_h\|_{L^2(\Omega)}.
	\end{align*}
	This, \eqref{ineq:G-best-approximation}--\eqref{ineq:inclusion}, the stability $\|\nabla \mathrm{J} \widetilde{e}_h\|_{L^2(\Omega)} \leq\|\mathrm{J}\|\|\widetilde e_h\|_{a,h}$ from
	\Cref{lem:enriching-operator}, a Cauchy inequality, and $(1-\alpha)\|\GrRec \widetilde{e}_h\|_{L^2(\Omega)}^2 + \beta\|h_\Tcal^{-1}S \widetilde{e}_h\|_{L^2(\Omega)}^{2}
	\leq
	\|\widetilde e_h\|_{a,h}^2$ from \eqref{def:ah} conclude the proof.
\end{proof}
\noindent The final lemma links \eqref{pr:discrete-problem} to \eqref{def:auxiliary-problem}. Recall the simple
eigenpair $(\lambda, u)$ of \eqref{pr:eigenvalue-problem} and the associated discrete eigenpair $(\lambda_h,u_h)$ of
\eqref{pr:discrete-problem} with $u_h = (u_\Tcal, u_\Fcal) \in V_h$ and $(u, u_\Tcal)_{L^2(\Omega)}\geq0$.
\begin{lemma}[upper bound for $\|u - u_\Tcal\|_{L^2(\Omega)}$]\label{lem:upper-bound-u-uT}
	If $h_{\max}$ is sufficiently small, then $\widetilde{u}_h = (\widetilde{u}_\Tcal, \widetilde{u}_\Fcal) \in V_h$ from \Cref{lem:best-approximation} with $f \coloneqq \lambda u$ satisfies
	\begin{align*}
	\|u - u_\Tcal\|_{L^2(\Omega)} \leq \cnst{cnst:spectral-gap} \|u - \widetilde{u}_\Tcal\|_{L^2(\Omega)}
	\end{align*}
	with the constant $\newcnst\label{cnst:spectral-gap} \coloneqq \sqrt{2}(1 + \max_{k \in \{1,\dots,N\}\setminus\{j\}} |\lambda/(\lambda_h(k) - \lambda)|) < \infty$.
\end{lemma}
\begin{proof}
	This follows as in \cite[Lem.\ 2.4]{CarstensenGallistlSchedensack2015} with straight-forward
	modifications and is hence omitted.
\end{proof}
\begin{proof}[Proof of \Cref{thm:L2-error-control}]
	The proof of \eqref{ineq:quasi-optimal-eigenvalue} is split into three steps.\medskip
	
	\noindent\emph{Step 1 provides the $L^2$ error estimate}
	\begin{align}
		\|u - u_\Tcal\|_{L^2(\Omega)} \leq \cnst{cnst:best-approximation-step-1-final} h_{\max}^s \min_{v_{p+1} \in P_{p+1}(\Tcal)}\|\nabla_\pw(u - v_{p+1})\|_{L^2(\Omega)}.
		\label{ineq:L2-estimate}
	\end{align}
	Recall $\widetilde{u}_h = (\widetilde{u}_\Tcal,\widetilde{u}_\Fcal)\in V_h$ from \Cref{lem:best-approximation} with
	$f \coloneqq \lambda u$.
	\Cref{lem:upper-bound-u-uT}, a triangle inequality, and \eqref{ineq:stability-2} with
	$\|(1-\Pi_{p+1})u\|_{L^2(\Omega)}\leq \|(1-G)u\|_{L^2(\Omega)}$ lead to
	\begin{align}
		\|u - u_\Tcal\|_{L^2(\Omega)} \leq \cnst{cnst:spectral-gap}\sigma_2 \|h_\Tcal (1 - \Pi_p) \nabla u\|_{L^2(\Omega)} +
		\cnst{cnst:spectral-gap}\|\Pi_{p+1} u - \widetilde{u}_\Tcal\|_{L^2(\Omega)}.
		\label{ineq:proof-L2-error-1}
	\end{align}
	Convergence rates for the error $\|\Pi_{p+1} u - \widetilde{u}_\Tcal\|_{L^2(\Omega)}$ in HHO methods for a source problem are established in \cite{DiPietroErnLemaire2014,DiPietroErn2015,ErnZanotti2020}.
	This proof follows \cite{CarstensenGallistlSchedensack2015,ErnZanotti2020} and utilizes the operator $\mathrm{J} :
	V_h \to V$ from \Cref{lem:enriching-operator}. Abbreviate $\widetilde{e}_h =
	(\widetilde{e}_\Tcal,\widetilde{e}_\Fcal) \coloneqq \I u - \widetilde{u}_h\in V_h$ and let $z \in V$ solve $-\Delta z = \widetilde{e}_\Tcal$ in $\Omega$, i.e., $z\in
	V$ satisfies
	\begin{align}
		(\nabla z, \nabla v)_{L^2(\Omega)} = (\widetilde{e}_\Tcal, v)_{L^2(\Omega)} \quad\text{for all } v \in V.
		\label{def:z}
	\end{align}
	Let $z_C \in S_0^{1}(\Tcal)\coloneqq P_{1}(\Tcal)\cap V$ denote the Scott-Zhang interpolation \cite{ScottZhang1990} of $z$
 	and observe that $(1-\Pi_p)\nabla z_C\equiv 0$ vanishes.
 	\Cref{lem:properties-reconstruction-operators} implies
	$S\I z_C\equiv
	0$ and therefore, the identity $a_h(\I z_C,\widetilde u_h) = (\nabla z_C, \GrRec \widetilde u_h)$ follows from
	\eqref{eq:ah-alternative-form} with $\GrRec \I z_C=\nabla z_C$.
	\Cref{lem:properties-reconstruction-operators} and $\I\mathrm J=1$ verify $\Pi_{\RT} \nabla \mathrm{J} \I u =
	\GrRec \I u=\Pi_\RT \nabla u$ and $\Pi_\RT \nabla \mathrm{J} \widetilde{u}_h = \GrRec \widetilde{u}_h$.
	This, $\nabla z_C\in P_0(\Tcal;\R^n) \subset \RT_p^{\pw}(\Tcal)$, and the symmetry of $a_h$ show
	\begin{align}\label{eqn:z_C}(\nabla z_C, \nabla \mathrm{J} \widetilde{e}_h)_{L^2(\Omega)}= (\nabla z_C, \nabla u - \GrRec
	\widetilde{u}_h)_{L^2(\Omega)} = a(u, z_C) - a_h(\widetilde u_h, \I z_C) = 0\end{align}
	with $a(u, z_C) = \lambda(u, z_C)_{L^2(\Omega)} = a_h(\widetilde u_h, \I z_C)$ from \eqref{pr:eigenvalue-problem} and
	\eqref{def:auxiliary-problem} in the last step.
	Hence, \eqref{def:z}--\eqref{eqn:z_C}, a Cauchy inequality, and
	$\|\nabla \mathrm{J} \widetilde e_h\|_{L^2(\Omega)} \leq\|\mathrm{J}\|\|\widetilde{e}_h\|_{a,h}$ from \Cref{lem:enriching-operator} confirm
	\begin{align}
		(\widetilde{e}_\Tcal, \mathrm{J} \widetilde{e}_h)_{L^2(\Omega)} =
		(\nabla (z-z_C), \nabla \mathrm{J} \widetilde{e}_h)_{L^2(\Omega)}
\leq\|\mathrm J\|\|\nabla(z - z_C)\|_{L^2(\Omega)} \|\widetilde{e}_h\|_{a,h}.
		\label{ineq:proof-L2-auxiliary}
	\end{align}
 	The stability estimate \eqref{ineq:stability-2} proves $\mathrm{osc}(\lambda u,\Tcal) \leq \lambda \sigma_2
	\|h_{\Tcal}^2(1 - \Pi_p) \nabla u\|_{L^2(\Omega)}$. This, \Cref{lem:best-approximation}, and \eqref{ineq:inclusion} provide
 	\begin{align}
		\|\widetilde{e}_h\|_{a,h} \leq \cnst{cnst:best-approximation-source}(1+\lambda\sigma_2h_{\max}^2) \|\nabla_\pw(u - Gu)\|_{L^2(\Omega)}.
 		\label{ineq:proof-L2-auxiliary-best-approximation}
 	\end{align}
	The elliptic regularity theory establishes $z \in V \cap H^{1+s}(\Omega)$ for $0 < s \leq 1$ on the polyhedral
	Lipschitz domain $\Omega$ and the approximation property of the Scott-Zhang interpolation $z_C$ \cite{ScottZhang1990}
	provides the constants $\newcnst\label{cnst:Scott-Zhang-approximation},\newcnst\label{cnst:elliptic-regularity}$ depending exclusively on the domain $\Omega$ such that
	\begin{align*}
		\cnst{cnst:Scott-Zhang-approximation}^{-1} h^{-s}_{\max}\|\nabla(z - z_C)\|_{L^2(\Omega)} \leq \|z\|_{H^{1+s}(\Omega)}
		\leq
		\cnst{cnst:elliptic-regularity} \|\widetilde{e}_\Tcal\|_{L^2(\Omega)}.
	\end{align*}
	Since $\Pi_{p+1}J\widetilde e_h=\widetilde e_\Tcal=\Pi_{p+1} u - \widetilde{u}_\Tcal$, the combination of \eqref{ineq:proof-L2-auxiliary}--\eqref{ineq:proof-L2-auxiliary-best-approximation}  
	verifies 
	\begin{align*}
		\|\Pi_{p+1} u - \widetilde{u}_\Tcal\|_{L^2(\Omega)}^2=(\widetilde e_\Tcal, J\widetilde e_h)_{L^2(\Omega)}
		\leq\cnst{cnst:best-approximation-step-1} h_{\max}^s \|\nabla_\pw(u
	- Gu)\|_{L^2(\Omega)}\|\widetilde e_\Tcal\|_{L^2(\Omega)}\end{align*}
	with $\newcnst\label{cnst:best-approximation-step-1}\coloneqq \|\mathrm{J}\|\cnst{cnst:best-approximation-source}\cnst{cnst:Scott-Zhang-approximation}\cnst{cnst:elliptic-regularity}(1+\lambda\sigma_2h_{\max}^2)$.
	This and \eqref{ineq:proof-L2-error-1}
	conclude the proof of \eqref{ineq:L2-estimate} with $\newcnst\label{cnst:best-approximation-step-1-final}\coloneqq \cnst{cnst:spectral-gap}(\sigma_2 h_{\rm
	max}^{1-s}+\cnst{cnst:best-approximation-step-1})$ and Step 1.\medskip

	\noindent\emph{Step 2 discusses the remaining term
	$
		|\lambda-\lambda_h|+ \|\I u-u_h\|^2_{a,h}$ on the left-hand side of \eqref{ineq:quasi-optimal-eigenvalue}}.
	Abbreviate $e_h \coloneqq \I u - u_h \in V_h$.
	Elementary algebra with the normalization $\|u\|_{L^2(\Omega)} = 1 = \|u_\Tcal\|_{L^2(\Omega)}$ reveals $2\lambda = \lambda\|u - u_\Tcal\|_{L^2(\Omega)}^2 + 2\lambda (u,u_\Tcal)_{L^2(\Omega)}$.
	This and $\|e_h\|_{a,h}^2 - \lambda_h = \|\I u\|_{a,h}^2 - 2a_h(\I u, u_h)$ result in
	\begin{align}
		&\lambda - \lambda_h + \|e_h\|_{a,h}^2\label{eq:proof-a-priori-eigenvalue-split}\\
		&\quad = \lambda\|u - u_\Tcal\|^2_{L^2(\Omega)} + \|\I u\|_{a,h}^2 - \lambda + 2(\lambda(u, u_\Tcal)_{L^2(\Omega)} - a_h(\I u, u_h)).\nonumber
	\end{align}
{{\emph{Step 2.1 bounds {$\|\I u\|_{a,h}^2 - \lambda$}.}}}
	The commutativity $\Pi_\RT \nabla u = \GrRec \I u$ from \Cref{lem:properties-reconstruction-operators} and 
	\eqref{eq:ah-alternative-form} with $\Pi_p \Pi_{\RT} \nabla u = \Pi_p \nabla u$ show \begin{align*}
		\|\I u\|_{a,h}^2 = (1 -
		\alpha)(\Pi_\RT \nabla u, \nabla u)_{L^2(\Omega)} + \alpha(\Pi_p \nabla u, \nabla u)_{L^2(\Omega)} + \linebreak\beta\|h_\Tcal^{-1} S \I u\|_{L^2(\Omega)}^2.
	\end{align*}
	This and $\lambda = \|\nabla u\|_{L^2(\Omega)}^2$ prove
	\begin{align*}
		\|\I u\|_{a,h}^2 - \lambda &= \|\I u\|_{a,h}^2 - \|\nabla u\|^2_{L^2(\Omega)} = -\alpha((1 - \Pi_p) \nabla u, \nabla u)_{L^2(\Omega)}\\
		&\quad - (1 - \alpha)((1 - \Pi_\RT)\nabla u, \nabla u)_{L^2(\Omega)} + \beta\|h_\Tcal^{-1} S \I u\|_{L^2(\Omega)}^2.
	\end{align*}
	Thus, $0<\alpha<1$ and %
	\eqref{ineq:stabilization-best-approximation} with
	$\widetilde{u}$ replaced by $u$ imply
	\begin{align}
		\|\I u\|_{a,h}^2 - \lambda \leq \beta\|h_\Tcal^{-1} S \I u\|_{L^2(\Omega)}^2\leq\beta C_P^2\|\nabla_\pw(u - Gu)\|_{L^2(\Omega)}^2.
		\label{ineq:proof-a-priori-eigenvalue-term-1}
	\end{align}
\emph{Step 2.2 controls $\lambda(u, u_\Tcal)_{L^2(\Omega)} - a_h(\I u, u_h)$.}
	The weak problem \eqref{pr:eigenvalue-problem} and $\Pi_{p+1}\mathrm J u_h=u_\Tcal$
	reveal
	\begin{align}
		&\lambda(u, u_\Tcal)_{L^2(\Omega)}
		\label{eq:proof-a-priori-eigenvalue-term-2-split}
		=  a(u, \mathrm{J} u_h) - \lambda ((1 - \Pi_{p+1}) u, \mathrm{J} u_h)_{L^2(\Omega)}.
	\end{align}
	\Cref{lem:properties-reconstruction-operators} provides $\Pi_\RT \nabla \mathrm{J} u_h = \GrRec u_h$ and $\GrRec \I
	u = \Pi_{\RT} \nabla u$. This and \eqref{eq:ah-alternative-form} lead to 
	\begin{align*}
		a_h(\I u, u_h) = ((1 -
		\alpha)\Pi_{\RT}\nabla u + \alpha \Pi_p \nabla u, \nabla \mathrm{J} u_h) + \beta(h_\Tcal^{-2} S \I u, S
		u_h)_{L^2(\Omega)}.
	\end{align*}
	This and \eqref{eq:proof-a-priori-eigenvalue-term-2-split} show
	\begin{align*}
		\lambda(u, u_\Tcal)_{L^2(\Omega)} &- a_h(\I u, u_h) =  - \lambda ((1 - \Pi_{p+1}) u, \mathrm{J} u_h)_{L^2(\Omega)}- \beta(h_\Tcal^{-2} S \I u, S
	u_h)_{L^2(\Omega)}\nonumber\\
		   &+(1-\alpha)((1 - \Pi_{\RT})\nabla u, \nabla \mathrm{J}u_h
		)_{L^2(\Omega)}
		+ \alpha((1 - \Pi_p)\nabla u, \nabla \mathrm J u_h)_{L^2(\Omega)}.
	\end{align*}
	Therefore, the Cauchy inequality and $P_p(\mathcal{T};\mathbb{R}^n) \subset \RT_p^\pw(\mathcal{T};\mathbb{R}^n)$ imply
	\begin{align}\label{eq:proof-a-priori-eigenvalue-term-2-split-2}
		\lambda(u, &u_\Tcal)_{L^2(\Omega)} - a_h(\I u, u_h) \leq \lambda \|(1 - \Pi_{p+1}) u\|_{L^2(\Omega)}\|(1 - \Pi_{p+1}) \mathrm{J} u_h\|_{L^2(\Omega)}\nonumber\\
		& + \|(1 - \Pi_p)\nabla u\|_{L^2(\Omega)}\|(1 - \Pi_{p})\nabla \mathrm J u_h\|_{L^2(\Omega)} - \beta(h_\Tcal^{-2} S \I u, S
		u_h)_{L^2(\Omega)}.
	\end{align}
	In the following, we control the terms on the right-hand side of \eqref{eq:proof-a-priori-eigenvalue-term-2-split-2}.
	The split $u_h=\I u-e_h$, $\|h_\Tcal^{-1} S\I u\|_{L^2(\Omega)}\geq0$, 
	and a Cauchy inequality provide
	\begin{align}\notag
		-(h_\Tcal^{-2} S \I u, S
		u_h)_{L^2(\Omega)}
		&\leq \|h_\Tcal^{-1}S\I u\|_{L^2(\Omega)}\|h_\Tcal^{-1} Se_h\|_{L^2(\Omega)}\\
		&\leq C_Pt/2\|\nabla_\pw(u-Gu)\|_{L^2(\Omega)}^2 + C_P/(2t)\|h_\Tcal^{-1}
		Se_h\|_{L^2(\Omega)}^2\label{eqn:S_term_bound_proof_5}
	\end{align}
	from \eqref{ineq:stabilization-best-approximation} with $\widetilde u$ replaced by $u$ and a Young inequality with
	arbitrary $t>0$ in the last step.
	Notice that $\Pi_p \nabla \mathrm{J}\I u = \mathcal{G} \I u = \Pi_p \nabla u$ by \Cref{lem:properties-reconstruction-operators}.
	Hence, a triangle inequality and $\|\nabla(u - \mathrm{J}\I u)\|_{L^2(\Omega)} \leq \cnst{cnst:enrichment-operator}\|\nabla_\pw(u - G u)\|_{L^2(\Omega)}$ from a combination of \eqref{ineq:J-best-approximation} with \eqref{ineq:G-best-approximation}--\eqref{ineq:inclusion} verify
	\begin{align}\label{ineq:proof-L2-estimate-1}
		\|(1-\Pi_p)\nabla \mathrm{J} \I
		u\|_{L^2(\Omega)} &\leq \|(1 - \Pi_p) \nabla u\|_{L^2(\Omega)} + \|\nabla(u - \mathrm{J}\I u)\|_{L^2(\Omega)}\nonumber\\
		&\leq (1+\cnst{cnst:enrichment-operator}) \|\nabla_\pw(u-Gu)\|_{L^2(\Omega)}.
	\end{align}
	This, \eqref{ineq:stability-2}, a triangle inequality with the split $u_h=\I u-e_h$, and the stability $\|\nabla \mathrm{J} e_h\|_{L^2(\Omega)} \leq\|\mathrm
	J\|\|e_h\|_{a,h}$ from \Cref{lem:enriching-operator} provide
	\begin{align}\label{ineq:proof-L2-estimate-2}
		\sigma_2^{-1}h_{\rm max}^{-1}&\|(1-\Pi_{p+1})\mathrm{J}u_h\|_{L^2(\Omega)} \leq \|(1 - \Pi_p)\nabla \mathrm{J} u_h\|_{L^2(\Omega)}\nonumber\\
	 	&\leq \|(1 - \Pi_p)\nabla \mathrm{J}\I u\|_{L^2(\Omega)} + \|(1 - \Pi_p)\nabla \mathrm{J} e_h\|_{L^2(\Omega)}\nonumber\\
	 	&\leq (1+\cnst{cnst:enrichment-operator}) \|\nabla_\pw(u-Gu)\|_{L^2(\Omega)} + \|\mathrm{J}\|\|e_h\|_{a,h}.
	\end{align}
	The combination of \eqref{ineq:proof-L2-estimate-1}--\eqref{ineq:proof-L2-estimate-2} with \eqref{ineq:stability-2} and a Young inequality with $t>0$ result in
	\begin{align}
		\notag\sigma_2^{-2}h_{\rm max}^{-2}&\|(1-\Pi_{p+1}) u\|_{L^2(\Omega)}\|(1-\Pi_{p+1})\mathrm{J}u_h\|_{L^2(\Omega)}\\&\leq
		\label{ineq:proof-a-priori-eigenvalue-term-2-J}
		\|(1-\Pi_{p})\nabla u\|_{L^2(\Omega)}\|(1-\Pi_{p})\nabla\mathrm{J}u_h\|_{L^2(\Omega)}\\&\leq
		(1+\cnst{cnst:enrichment-operator}+\|\mathrm{J}\|t/2)\|\nabla_\pw(u-Gu)\|_{L^2(\Omega)}^2+\|\mathrm{J}\|/(2t)\|e_h\|_{a,h}^2.\notag
	\end{align}
	Then \eqref{eq:proof-a-priori-eigenvalue-term-2-split-2}--\eqref{eqn:S_term_bound_proof_5}, \eqref{ineq:proof-L2-estimate-2}--\eqref{ineq:proof-a-priori-eigenvalue-term-2-J} with the choice
	$t\coloneqq 2(1+\lambda\sigma_2^2h_{\rm max}^2)\|\mathrm{J}\|+2C_P$, and
	$\beta\|h_\Tcal^{-1}Se_h\|_{L^2(\Omega)}^2\leq\|e_h\|_{a,h}^2$ from \eqref{def:ah} lead to
	\begin{align}\label{eqn:last_5.3}
		\lambda(u, u_\Tcal)_{L^2(\Omega)} - a_h(\I u, u_h) - \|e_h\|_{a,h}^2/4 \leq \cnst{cnst:best-approximation-step-2-2} \|\nabla_\pw(u - Gu)\|_{L^2(\Omega)}^2
	\end{align}
	with $\newcnst\label{cnst:best-approximation-step-2-2}\coloneqq (1+\lambda\sigma_2^2h_{\rm max}^2)(1+\cnst{cnst:enrichment-operator} + \|\mathrm{J}\|t/2) + \beta C_Pt/2$.
	\medskip

	\noindent\emph{Step 3 finishes the proof.}
	\Cref{thm:leb} guarantees $\lambda_h \leq \lambda$ for sufficiently small mesh-sizes $h_{\max} \leq (\alpha/(\lambda \sigma_2^2))^{1/2}$.
	This, the combination of 
	\eqref{eq:proof-a-priori-eigenvalue-split}--\eqref{ineq:proof-a-priori-eigenvalue-term-1} and \eqref{eqn:last_5.3}
	with the $L^2$ error
	estimate \eqref{ineq:L2-estimate} from Step 1 results in 
	\begin{align*}
		|\lambda-\lambda_h|+ \|e_h\|_{a,h}^2/2\leq(\lambda\cnst{cnst:best-approximation-step-1-final}^{{2}} h_{\rm max}^{2s} +\beta C_P^2 + 2\cnst{cnst:best-approximation-step-2-2}) \|\nabla_\pw(u-Gu)\|_{L^2(\Omega)}^2.
	\end{align*}
Thus, \eqref{ineq:L2-estimate} and \eqref{ineq:G-best-approximation} conclude the proof of
\eqref{ineq:quasi-optimal-eigenvalue} with $\cnst{cnst:best-approximation} \coloneqq 2(\lambda\cnst{cnst:best-approximation-step-1-final}^{{2}} h_{\rm max}^{2s}
+\beta C_P^2 + 2\cnst{cnst:best-approximation-step-2-2})+\cnst{cnst:best-approximation-step-1-final}^2$.
\end{proof}
\Cref{thm:L2-error-control} implies the following convergence rates and recovers \cite{CaloCicuttinDengErn2019,CarstensenErnPuttkammer2021} for the eigenvalues and eigenfunctions error in the $H^1$ seminorm.
\begin{corollary}[convergence]
	If $u \in V \cap H^{1+t}(\Omega)$ for $s \leq t \leq p+1$, then 
	\begin{align*}
		h_{\max}^{-s}\|u - u_\Tcal\|_{L^2(\Omega)}+h_{\max}^{-t}\left(|\lambda - \lambda_h| + \|\I u -
		u_h\|_{a,h}^2\right) \lesssim h_{\max}^{t}~\text{as }h_{\max}\to 0.
	\end{align*}
\end{corollary}
\begin{proof}
	This follows immediately from \Cref{thm:L2-error-control}, the stability \eqref{ineq:stability-1}, and standard
	approximation properties of piecewise polynomials \cite[Lemma 4.3.8]{BrennerScott2008}.
\end{proof}
The techniques of this section also apply to the HHO method of \cite{CarstensenErnPuttkammer2021} and lead to the
optimal rate $s + t$ for the $L^2$ error towards a simple eigenvalue therein.
\section{A posteriori error analysis}\label{sec:a-posteriori}
The two assumptions \ref{def:A1}--\ref{def:A2} below concern some $q \in
H^1(\Tcal;\R^n)$ and lead to a stabilization-free a~posteriori error control of $\|\nabla u - q\|_{L^2(\Omega)}$ in two
or three space dimensions. Let $\RT_0(\Tcal)\coloneqq \RT_0^\pw(\Tcal)\cap
H(\div)$ denote the lowest-order conforming Raviart-Thomas space, set $S^m_0(\Tcal)\coloneqq P_m(\Tcal)\cap
H^1_0(\Omega)$ for $m\in\mathbb N$, and
suppose
\begin{enumerate}[wide,label = (A\arabic*)]
	\item\label{def:A1} $(q, \nabla v_C)_{L^2(\Omega)} = \lambda_h(u_\Tcal, v_C)_{L^2(\Omega)}$ for all $v_C \in S_{0}^{1}(\Tcal)$,
	\item\label{def:A2} $(q, q_\RT)_{L^2(\Omega)} = 0$ for all $q_\RT \in \RT_0(\Tcal)$ with $\div\, q_\RT = 0$.
\end{enumerate}

\begin{theorem}[a~posteriori]\label{thm:a-posteriori}
	Any $q \in H^1(\Tcal;\R^n)$ with \ref{def:A1}--\ref{def:A2} and $\eta$ from \eqref{def:eta-local} with $p_h$ replaced by $q$ satisfy
	\begin{align}
		\cnst{cnst:a-posteriori}^{-1}\|\nabla u - q\|_{L^2(\Omega)}^2 \leq \eta^2 + \|\lambda u - \lambda_h u_\Tcal\|^2_{L^2(\Omega)}.
		\label{ineq:a-posteriori-reliability}
	\end{align}
	The constant $\newcnst\label{cnst:a-posteriori}$ only depends on $p$, $n$, $\Omega$, and the shape regularity of $\mathcal{T}$.
\end{theorem}
\begin{proof}
	This is an extension of \cite{BertrandCarstensenGraessle2021} to eigenvalue problems. For the convenience of the reader, the main arguments are briefly outlined below.
	Let $\psi \in V$ solve $- \Delta \psi = -\div\,q \in H^{-1}(\Omega)$ so that the Pythagoras theorem allows for the split
	\begin{align}
		\|\nabla u - q\|_{L^2(\Omega)}^2 = \|\nabla(u - \psi)\|_{L^2(\Omega)}^2 + \|\nabla \psi - q\|_{L^2(\Omega)}^2.
		\label{eq:a-posteriori-BCG-split}
	\end{align}
	\emph{Upper bound for $\|\nabla(u - \psi)\|_{L^2(\Omega)}$.}
	Abbreviate $\varrho \coloneqq u - \psi \in V$ and let $\varrho_C \in S^1_0(\Tcal)$ denote the Scott-Zhang interpolation of $\varrho$ \cite{ScottZhang1990}. Then \ref{def:A1}, $(\nabla \psi, \nabla \varrho)_{L^2(\Omega)} = (q, \nabla \varrho)_{L^2(\Omega)}$, and \eqref{pr:eigenvalue-problem} lead to
	\begin{align}
		\|\nabla \varrho\|^2_{L^2(\Omega)} &= \lambda b(u, \varrho) - \lambda_h(u_\Tcal, \varrho_C)_{L^2(\Omega)} - (q,\nabla (\varrho - \varrho_C))_{L^2(\Omega)}\nonumber\\
		&= (\lambda u - \lambda_h u_\Tcal, \varrho)_{L^2(\Omega)} + \lambda_h(u_\Tcal, \varrho - \varrho_C)_{L^2(\Omega)} - (q,\nabla(\varrho - \varrho_C))_{L^2(\Omega)}.\label{eq:a-posteriori-term-1-integration-by-parts}
	\end{align}
	The last two $L^2$ scalar products on the right-hand side of \eqref{eq:a-posteriori-term-1-integration-by-parts}
	arise in the explicit residual-based a posteriori error estimation of standard conforming FEM for the Poisson model problem, cf., e.g., \cite[Section 2.2]{AinsworthOden2000} or \cite[Chapter 34]{ErnGuermond2021}, and are controlled by
	\begin{align*}
		\Big(\|h_\Tcal(\div_\pw q + \lambda_h u_\Tcal)\|_{L^2(\Omega)}^2 + \sum_{F \in \Fcal(\Omega)} h_F\|[q\cdot \nu_F]_F\|_{L^2(F)}^2\Big)^{{1/2}}\|\nabla \varrho\|_{L^2(\Omega)}.
	\end{align*}
	This, \eqref{eq:a-posteriori-term-1-integration-by-parts}, a Cauchy inequality, and a Friedrichs inequality result in
	\begin{align}
		\|\nabla(u - \psi)\|_{L^2(\Omega)}^2 \lesssim \eta^2 + \|\lambda u - \lambda_h u_\Tcal\|^2_{L^2(\Omega)}.
		\label{ineq:a-posteriori-BCG-term-1}
	\end{align}
	\emph{Upper bound for $\|\nabla \psi - q\|_{L^2(\Omega)}$.}
	The function $\phi \coloneqq \nabla \psi - q \in L^2(\Omega;\R^n)$ is divergence-free $\div\,\phi = 0$ and
	orthogonal to the divergence-free Raviart-Thomas functions $q_{\RT}\in\RT_0(\Tcal)$ from \ref{def:A2}.
	The Helmholtz decomposition on a simply connected domain $\Omega$ immediately implies
	$\mathrm{Curl}\,\beta = \phi$ for some $\beta\in H^1(\Omega;\R^{2n-3})$, but
in this paper, the domain $\Omega$ does not need to be simply connected.
	However, the extra condition \ref{def:A2} ensures the existence 
	of some orthogonal correction $\phi_{\RT}\in\RT_0(\Tcal)$ with $\div\,\phi_{\RT}=0$
	such that
	the integrals $\int_{\Gamma_j} (\phi-\phi_{\RT})\cdot\nu\d s=0$ over the
	$J\in\mathbb N$ connectivity components $\Gamma_j$ for $j=1,...,J$ of $\partial \Omega$ vanish, cf.~\cite[Lemma
		2]{BertrandCarstensenGraessle2021} for further details.
	Thus classical theorems \cite{GiraultRaviart1986} imply the existence of $\beta\in H^1(\Omega;\R^{2n-3})$ such that
	$\mathrm{Curl}\,\beta = \phi-\phi_{\RT}$ and
	$\|\nabla \beta\|_{L^2(\Omega)}\lesssim \|\phi\|_{L^2(\Omega)}$.
	Since the Scott-Zhang interpolation $\beta_C \in S^1_0(\Tcal;\R^{2n-3})$ of $\beta$ satisfies $\mathrm{Curl}\,\beta_C \in \RT_0(\Tcal)$ and $\div\, \mathrm{Curl}\,\beta_C = 0$, \ref{def:A2} shows
	\begin{align*}
		\|\nabla \psi - q\|_{L^2(\Omega)}^2 = (\phi,\mathrm{Curl}\,\beta+\phi_{\RT})_{L^2(\Omega)}= (\phi,\mathrm{Curl}(\beta - \beta_C))_{L^2(\Omega)}.
	\end{align*}
	A piecewise integration by parts, the trace inequality, the approximation property of the Scott-Zhang interpolation \cite{ScottZhang1990}, and the Cauchy inequality lead to
	\begin{align}
		\|\nabla \psi - q\|_{L^2(\Omega)}^2 \lesssim \|h_\Tcal \mathrm{curl}\,q\|_{L^2(\Omega)}^2 + \sum_{F \in \Fcal}
		h_F\|[q \times \nu_F]_F\|_{L^2(F)}^2.
		\label{ineq:a-posteriori-BCG-term-2}
	\end{align}
	The combination of \eqref{eq:a-posteriori-BCG-split} with \eqref{ineq:a-posteriori-BCG-term-1}--\eqref{ineq:a-posteriori-BCG-term-2} concludes the proof of \eqref{ineq:a-posteriori-reliability}.
\end{proof}
\noindent One key observation is that $q \coloneqq p_h \coloneqq \Pi_p \GrRec u_h$ satisfies \ref{def:A1}--\ref{def:A2} {as shown in the proof of \Cref{thm:rel-eff} below}.
This leads to reliable a~posteriori error control for $\|\nabla u - p_h\|_{L^2(\Omega)}$. \Cref{thm:a-posteriori} can also be applied to the HHO scheme of \cite{CaloCicuttinDengErn2019}, where $q \coloneqq \nabla_\pw \PotRec u_h$ satisfies \ref{def:A1}--\ref{def:A2} for $p \geq 1$.
The lowest-order case $p=0$ therein can be treated separately as in \cite{BertrandCarstensenGraessle2021}.
\begin{theorem}[reliability and efficiency]\label{thm:rel-eff}
	For sufficiently small mesh-sizes $h_{\max}$, $p_h \coloneqq \Pi_p \GrRec u_h \in P_p(\Tcal;\R^n)$ and $\eta$ from \eqref{def:eta-local} satisfy \eqref{ineq:a-posteriori}.
	The constants $C_\mathrm{eff}$ and $C_\mathrm{rel}$ exclusively depend on $p$, $n$, $\Omega$, and the shape regularity of $\mathcal{T}$.
\end{theorem}
\begin{proof}
	The first part of the proof verifies that $p_h = \Pi_p \GrRec u_h$ satisfies \ref{def:A1}--\ref{def:A2}.\medskip
	
	\noindent\emph{Proof of \ref{def:A1}.}
	Any $v_C \in S^1_0(\Tcal)$ satisfies $\nabla v_C = \GrRec \I v_C \in P_0(\Tcal)$ and $v_C = \PotRec \I v_C$. Thus $S \I v_C = 0$ and so,
	\begin{align*}
		(p_h,\nabla v_C)_{L^2(\Omega)} = (\GrRec u_h,\nabla v_C)_{L^2(\Omega)} = a_h(u_h,\I v_C) = \lambda_h(u_\Tcal,v_C)_{L^2(\Omega)}.
	\end{align*}
	\emph{Proof of \ref{def:A2}.} Given $q_\RT \in \RT_0(\Tcal) \subset H(\div,\Omega)$ with $\div\, q_\RT = 0$, the
	normal jump $[q_\RT \cdot \nu_F]_F$ vanishes on any interior side $F \in \Fcal(\Omega)$. Since divergence-free
	functions in $\RT_0(\Tcal)$ are piecewise constant, the definition of $\GrRec$ from
	\eqref{def:gradient-reconstruction} shows $(p_h,q_\RT)_{L^2(\Omega)} = (\GrRec u_h,q_\RT)_{L^2(\Omega)} = 0$ and
	concludes the proof of \ref{def:A2}.\medskip
	
	\noindent\emph{Proof of reliability.}
	Since $q = p_h$ satisfies \ref{def:A1}--\ref{def:A2}, \Cref{thm:a-posteriori} asserts
	\begin{align}
		\cnst{cnst:a-posteriori}^{-1}\|\nabla u - p_h\|_{L^2(\Omega)}^2 \leq \eta^2 +  \|\lambda u - \lambda_h u_\Tcal\|^2_{L^2(\Omega)}.
		\label{ineq:proof-a-posteriori-reliability}
	\end{align}
	The normalization $\|u\|_{L^2(\Omega)} = 1 = \|u_\Tcal\|_{L^2(\Omega)}$, elementary algebra, and the combination of
	the a~priori estimate
	\eqref{ineq:quasi-optimal-eigenvalue} with \eqref{ineq:G-best-approximation} reveal
	\begin{align}
		\|\lambda u - \lambda_h u_\Tcal\|^2_{L^2(\Omega)}
		& = (\lambda - \lambda_h)^2 + \lambda\lambda_h\|u -
		u_\Tcal\|^2_{L^2(\Omega)}\nonumber\\
		& \leq \cnst{cnst:a-posteriori-absorption}h_{\max}^{2s} \|\nabla_\pw(u  - G u)\|_{L^2(\Omega)}^2
		\label{ineq:proof-a-posteriori-superconvergence}
	\end{align}
	with the elliptic regularity of $u \in V \cap H^{1+s}(\Omega)$ for the parameter $0 < s \leq 1$ and 
	$\newcnst\label{cnst:a-posteriori-absorption}\coloneqq\max\{|\lambda-\lambda_h|, \lambda\lambda_h\}\cnst{cnst:best-approximation}$.
	The inequalities \eqref{ineq:stability-1} and $p_h\in P_p(\Tcal; \R^n)$ prove
	\begin{align}
		C_\mathrm{st,2}^{-1}\|\nabla_\pw(u  - G u)\|_{L^2(\Omega)} \leq \|(1 - \Pi_p) \nabla u\|_{L^2(\Omega)}
		\leq \|\nabla u - p_h\|_{L^2(\Omega)}.
		\label{ineq:R=GI-best-approximation}
	\end{align}
	For sufficiently small mesh-sizes $h_{\max}$, $\newcnst\label{cnst:a-posteriori-absorption-small-mesh-size}\coloneqq \cnst{cnst:a-posteriori}\cnst{cnst:a-posteriori-absorption}h_{\rm
	max}^{2s}C_{\rm st,2}^2<1$ and \eqref{ineq:proof-a-posteriori-reliability}--\eqref{ineq:R=GI-best-approximation} lead to
	\begin{align}
		\|\nabla u - p_h\|_{L^2(\Omega)}^2 \leq \cnst{cnst:a-posteriori} (1-\cnst{cnst:a-posteriori-absorption-small-mesh-size})^{-1} \eta^2.
		\label{ineq:reliability}
	\end{align}
	Under the additional assumption $h_{\max} \leq (\alpha/(\lambda \sigma_2^2))^{1/2}$, the quasi-best approximation \eqref{ineq:quasi-optimal-eigenvalue} and \eqref{ineq:R=GI-best-approximation}--\eqref{ineq:reliability} conclude the proof of
	\begin{align}
		|\lambda - \lambda_h| + \|\I u - u_h\|_{a,h}^2 + \|\nabla u - p_h\|_{L^2(\Omega)}^2 \leq C_{\rm rel} \eta^2
		\label{ineq:reliability-eigenvalue}
	\end{align}
	with $C_{\rm rel}\coloneqq (1+\cnst{cnst:best-approximation}C_{\rm st, 2}^2)\cnst{cnst:a-posteriori} (1-\cnst{cnst:a-posteriori-absorption-small-mesh-size})^{-1}$.\medskip
	
	\noindent\emph{Proof of efficiency.} The proof of $\eta^2 \lesssim \|\nabla u - p_h\|_{L^2(\Omega)}^2$ utilizes bubble-function techniques from \cite{Verfuerth1986}. Similar arguments are employed in \cite{DariDuranPadra2012} for the Crouzeix-Raviart FEM and, e.g., in \cite{AinsworthOden2000,DiPietroTittarelli2018,BertrandCarstensenGraessle2021,ErnGuermond2021} for the source problem. The efficiency $\sum_{F \in \Fcal} h_F\|[p_h \times \nu_F]_F\|_{L^2(F)}^2 \lesssim \|\nabla u - p_h\|^2_{L^2(\Omega)}$ 
	follows from the arguments in the proof of \cite[Lemma 7]{BertrandCarstensenGraessle2021} for the
	Poisson model problem; hence further details are omitted. The focus is therefore on the proof of the efficiency of
	\begin{align*}
		\|h_\Tcal \mathrm{curl}\,p_h\|_{L^2(\Omega)}^2 + \|h_\Tcal(\div\, p_h + \lambda_h u_\Tcal)\|^2_{L^2(\Omega)} + \sum_{F \in \Fcal(\Omega)} h_F\|[p_h \cdot \nu_F]_F\|^2_{L^2(F)}.
	\end{align*}
	Given $F \in \Fcal(\Omega)$, let $b_F \in S^n(\Tcal)$ denote the face-bubble function with $0 \leq b_F \leq 1$ in
	$\Omega$ and $\mathrm{supp}(b_F) = \overline{\omega_F}$ \cite[Section 3.1]{Verfuerth1986}. Define $\varrho \in
	S^{p+n}_0(\Tcal)$ such that $\varrho|_F = b_F[p_h\cdot \nu_F]_F \in P_{p+n}(F)$, $\mathrm{supp}(\varrho) =
	\overline{\omega_F}$, and $\varrho$ vanishes at all Lagrange points \cite{BrennerScott2008} in $\overline{\Omega}\setminus F$. Inverse estimates \cite[Ineq.~(3.2)]{Verfuerth1986} and an integration by parts prove, for any $F \in \Fcal(\Omega)$, that
	\begin{align*}
		\|[p_h \cdot \nu_F]_F\|_{L^2(F)}^2 \lesssim (\varrho, [p_h \cdot \nu_F]_F)_{L^2(F)} = (\nabla \varrho, p_h)_{L^2(\omega_F)} + (\varrho, \div_\pw p_h)_{L^2(\omega_F)}.
	\end{align*}
	This, $(\nabla u, \nabla \varrho)_{L^2(\omega_F)} = \lambda(u,\varrho)_{L^2(\Omega)},$ and a Cauchy inequality imply
	\begin{align*}
		&\|[p_h \cdot \nu_F]_F\|_{L^2(F)}^2 \lesssim \|\nabla \varrho\|_{L^2(\omega_F)}\|\nabla u - p_h\|_{L^2(\omega_F)}\\
		&\qquad + \|\varrho\|_{L^2(\omega_F)}\|\lambda u - \lambda_h u_\Tcal\|_{L^2(\omega_F)} + \|\varrho\|_{L^2(\omega_F)}\|\div_\pw p_h + \lambda_h u_\Tcal\|_{L^2(\omega_F)}.
	\end{align*}
	This, the inverse estimate $\|\nabla \varrho\|_{L^2(\omega_F)} \lesssim h_F^{-1}\|\varrho\|_{L^2(\omega_F)}$ \cite[Lemma 4.5.3]{BrennerScott2008}, and $\|\varrho\|_{L^2(\omega_F)}^2 \approx h_F\|\varrho\|_{L^2(F)}^2 \leq h_F\|[p_h\cdot \nu_F]_F\|_{L^2(F)}^2$ \cite[Ineq.~(3.5)]{Verfuerth1986} show
	\begin{align}
		&h_F\|[p_h \cdot \nu_F]_F\|^2_{L^2(F)} \lesssim \|\nabla u - p_h\|_{L^2(\omega_F)}^2\label{ineq:efficiency-normal-jump}\\
		&\qquad + h_F^2\|\lambda u - \lambda_h u_\Tcal\|_{L^2(\omega_F)}^2 + h_F^2\|\div_\pw p_h + \lambda_h u_\Tcal\|_{L^2(\omega_F)}^2.\nonumber
	\end{align}
	Let $b_T \in P_{n+1}(T) \cap W^{1,\infty}_0(T)$ denote the volume-bubble function in $T\in\Tcal$ with $0 \leq b_T \leq 1$ and $b_T = 0$ on $\partial T$ \cite[Section 3.1]{Verfuerth1986}.
	Abbreviate $v_{p+1} \coloneqq \div\, p_h + \lambda_h u_T \in P_{p+1}(T)$ and
	define $\varphi\coloneqq b_T v_{p+1} \in S^{p+n+2}_0(T)\coloneqq P_{p+n+2}(T)\cap H^1_0(T)\subset V$.
	Inverse estimates \cite[Ineq.~(3.1)]{Verfuerth1986} and an integration by parts imply
	\begin{align}
		\|v_{p+1}\|_{L^2(T)}^2 \lesssim (\varphi, v_{p+1})_{L^2(T)} = -(\nabla \varphi, p_h)_{L^2(T)} + (\varphi, \lambda_h u_T)_{L^2(T)}.
		\label{ineq:efficiency-volume-1}
	\end{align}
	Since (the extension by zero of) $\varphi$ belongs to $V$, \eqref{pr:eigenvalue-problem} provides $(\nabla \varphi, \nabla u)_{L^2(T)} =
	\lambda(u, \varphi)_{L^2(T)}$. This,
	\eqref{ineq:efficiency-volume-1}, and a Cauchy inequality lead to
	\begin{align*}
		\|v_{p+1}\|_{L^2(T)}^2 \lesssim \|\nabla \varphi\|_{L^2(T)}\|\nabla u - p_h\|_{L^2(T)} + \|\varphi\|_{L^2(T)}\|\lambda u - \lambda_h u_T\|_{L^2(T)}.
	\end{align*}
	Hence $\|\varphi\|_{L^2(T)} = \|b_T v_{p+1}\|_{L^2(T)} \leq \|v_{p+1}\|_{L^2(T)}$ from $0 \leq b_T \leq 1$ in $T$ and the inverse estimate
	$\|\nabla \varphi\|_{L^2(T)} \leq h_T^{-1}\|\varphi\|_{L^2(T)}$ \cite[Lemma 4.5.3]{BrennerScott2008} reveal
	\begin{align}
		h_T^2\|\div\, p_h + \lambda_h u_T\|^2_{L^2(T)} \lesssim \|\nabla u - p_h\|_{L^2(T)}^2 + h_T^2\|\lambda u - \lambda_h u_T\|_{L^2(T)}^2.
		\label{ineq:efficiency-volume-div}
	\end{align}
	The local estimate $h_T\|\mathrm{curl}\,p_h\|_{L^2(T)} \lesssim \|\nabla u - p_h\|_{L^2(T)}$ follows from similar
	arguments as above and details are omitted.
	The combination of this with the local estimates \eqref{ineq:efficiency-normal-jump} and \eqref{ineq:efficiency-volume-div} results in
	$\eta^2 \lesssim \|\nabla u - p_h\|_{L^2(\Omega)}^2 + \|h_\Tcal(\lambda u - \lambda_h u_\Tcal)\|_{L^2(\Omega)}^2$.
	This and the control over $\|\lambda u - \lambda_h u_\Tcal\|_{L^2(\Omega)}$ in
	\eqref{ineq:proof-a-posteriori-superconvergence}--\eqref{ineq:R=GI-best-approximation} lead to the efficiency $\eta^2 \lesssim \|\nabla u - p_h\|_{L^2(\Omega)}^2$.
\end{proof}
	
\section{Numerical examples}\label{sec:numerical-examples}
The section presents three numerical benchmarks for the approximation of Dirichlet eigenvalues of the Laplacian on nonconvex domains $\Omega \subset \R^2$.
\subsection{Parameter selection}
\label{sec:Parameter selection}
For right-isosceles triangles, recall $C_\mathrm{st,2} \leq \sqrt{2}$ from \Cref{ex:C_st} and $C_P = 1/(\sqrt{2}\pi)$
from \cite{LaugesenSiudeja2010}. Throughout this section, let $\alpha = 0.5$ and $\beta \coloneqq \alpha/\sigma_2^2 =
4.934802$ with $\sigma_2^2 = C_P^2 C_{\mathrm{st,2}}^2 = 1/\pi^2 = 0.101321$.
The computable (a~posteriori) condition $\sigma_2^2 \max\{\beta, h_{\max}^2 \lambda_h(j)\} \leq \alpha$ from \Cref{thm:leb} leads to $\mathrm{GLB}(j) \coloneqq \lambda_h(j) \leq \lambda(j)$.
Since the parameters are chosen before-hand, the condition $h_{\max}^2 \lambda_h \leq \alpha/\sigma_2^2 = 4.934802$ may
not be satisfied on a coarse mesh with large $h_{\max}$ and $j$. In this case, $\mathrm{GLB}(j) \coloneqq 0$ (which is a guaranteed lower eigenvalue bound), so only GLB are displayed in this {section}.

\subsection{Numerical realization}
The algebraic eigenvalue problem \eqref{eq:algebraic-eigenvalue-problem} is realized with the iterative solver
\texttt{eigs} from the MATLAB standard library in an extension of the data structures and short MATLAB programs in
\cite{AlbertyCFunken1999,CBrenner2017}; the termination and round-off errors are expected to be very small and neglected
for simplicity.\par
\begin{figure}
	\centering
	\begin{minipage}[t]{0.33\textwidth}
		\centering
		\includegraphics[scale=0.5]{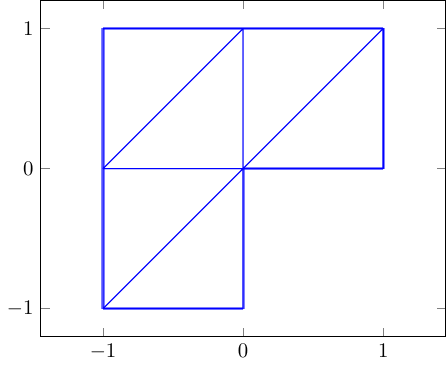}
	\end{minipage}\hfill
	\begin{minipage}[t]{0.33\textwidth}
		\centering
		\includegraphics[scale=0.5]{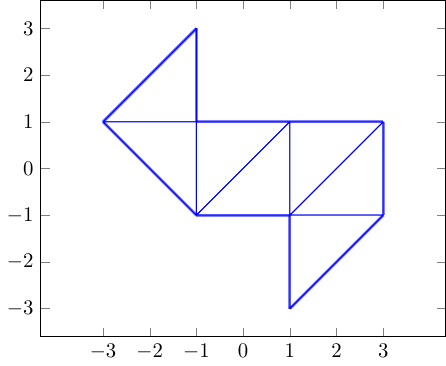}
	\end{minipage}\hfill
	\begin{minipage}[t]{0.33\textwidth}
		\centering
		\includegraphics[scale=0.5]{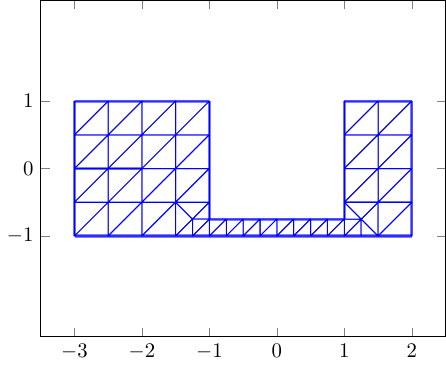}
	\end{minipage}\hfill
	\caption{Initial triangulations of the L-shaped domain in \Cref{sec:L-shape}, the isospectral drum in
		\Cref{sec:isospectral-domain}, and the dumbbell-slit domain in \Cref{sec:dumbbellslit}.}
	\label{fig:initial-triangulation}
\end{figure}
\begin{figure}
	\centering
	\includegraphics{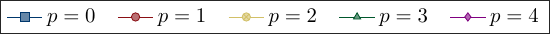}
	\caption{Polynomial degrees $p = 0, \dots, 4$ in the numerical benchmarks of \Cref{sec:numerical-examples}.}
	\label{fig:legend}
\end{figure}
The a~posteriori estimate from \Cref{thm:a-posteriori} motivates the refinement indicator $\eta^2(T)$ from
\eqref{def:eta-local} with $\eta^2 = \sum_{T \in \Tcal} \eta^2(T)$. The standard adaptive algorithm \cite[Algorithm 2.2]{AFEM} is modified in that, if $h_{\max}^2 \lambda_h \leq \alpha/\sigma_2^2$ is not satisfied, the mesh is uniformly refined.
It runs with the initial triangulations from \Cref{fig:initial-triangulation}, the default bulk parameter $\theta = 0.5$, and polynomial degrees $p$ displayed in \Cref{fig:legend}. 
\par
The uniform and adaptive mesh-refinements lead to convergence history plots of the eigenvalue error $\lambda(j) - \mathrm{GLB}(j)$ and the a~posteriori estimate $\eta^2$ plotted against the number of degrees of freedom of $V_h$ (ndof) in log-log plots below; dashed (resp.~solid) lines represent uniform (resp.~adaptive) mesh-refinements.
\subsection{L-shaped domain}\label{sec:L-shape}
\begin{figure}[]
	\centering
	\begin{minipage}[t]{0.49\textwidth}
		\centering
		\includegraphics[scale=0.8]{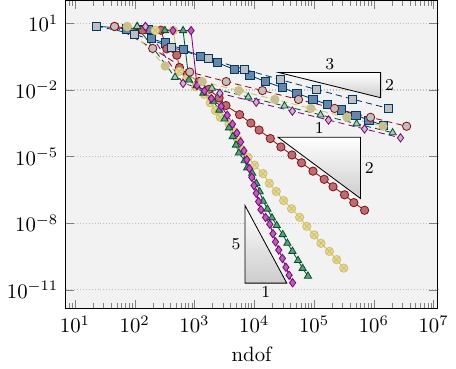}
	\end{minipage}\hfill
	\begin{minipage}[t]{0.49\textwidth}
		\centering
		\includegraphics[scale=0.8]{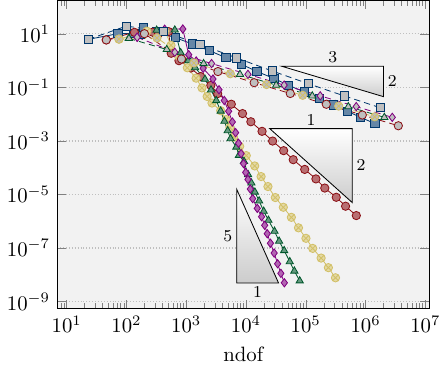}
	\end{minipage}\hfill
	\caption{Convergence history plot of $\lambda(1) - \mathrm{GLB}(1)$ (left) and $\eta^2$ (right) for polynomial
		degrees $p$ from \Cref{fig:legend} on uniform (dashed line) and adaptive (solid line) triangulations of the L-shaped
		domain in \Cref{sec:L-shape}.}
	\label{fig:L-shape}
\end{figure}

\begin{figure}
	\centering
	\begin{minipage}[t]{0.49\textwidth}
		\centering
		\includegraphics[scale=0.8]{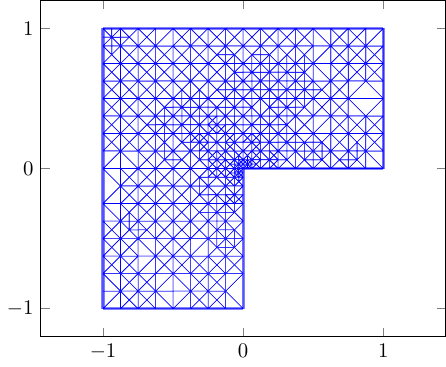}
	\end{minipage}\hfill
	\begin{minipage}[t]{0.49\textwidth}
		\centering
		\includegraphics[scale=0.8]{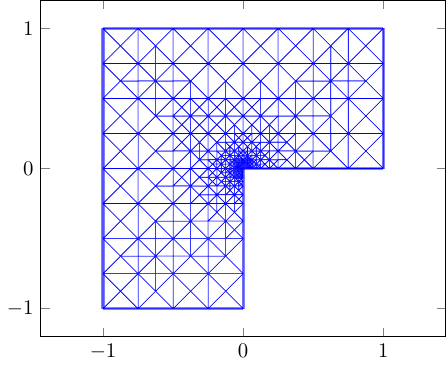}
	\end{minipage}\hfill
	\caption{Adaptive triangulations of the L-shaped domain in \Cref{sec:L-shape} into 1034 triangles for $p = 0$ (left)
	and into 1138 triangles for $p = 3$ (right) for the approximation to $\lambda(1)$.}
	\label{fig:L-shape-triangulation}
\end{figure}

The first example concerns the principle Dirichlet eigenvalue on the domain $\Omega \coloneqq (-1,1)^2\setminus ([0,1)\times[0,-1))$ with a re-entering corner at $(0,0)$ and the reference value $\lambda(1) = 9.6397238440219410$ from \cite{BetckeTrefethen2005}. 
This leads to the suboptimal convergence rate $2/3$ for $\lambda(1) - \mathrm{GLB}(1)$ and $\eta^2$ (for all $p$) on uniform triangulations in \Cref{fig:L-shape}.
The adaptive mesh-refining algorithm refines towards the origin as displayed in \Cref{fig:L-shape-triangulation} and recovers the optimal convergence rates $p+1$ for $\lambda(1) - \mathrm{GLB}(1)$ and $\eta^2$.

\subsection{Isospectral domain}\label{sec:isospectral-domain}
\begin{figure}
	\centering
	\begin{minipage}[t]{0.49\textwidth}
		\centering
		\includegraphics[scale=0.8]{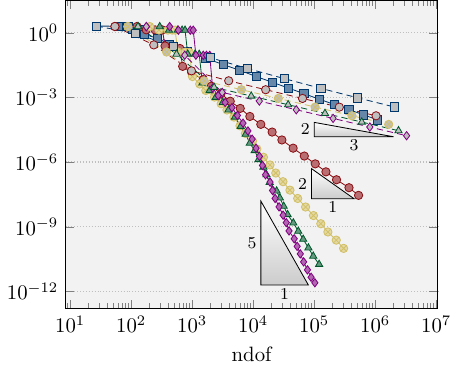}
	\end{minipage}\hfill
	\begin{minipage}[t]{0.49\textwidth}
		\centering
		\includegraphics[scale=0.8]{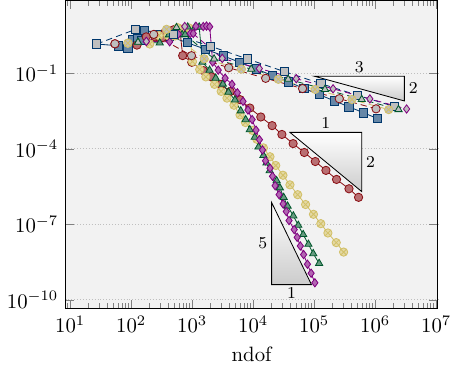}
	\end{minipage}\hfill
	\caption{Convergence history plot of $\lambda(1) - \mathrm{GLB}(1)$ (left) and $\eta^2$ (right) for polynomial
	degrees $p$ from \Cref{fig:legend} on uniform (dashed line) and adaptive (solid line) triangulations of the
isospectral domain in \Cref{sec:isospectral-domain}.}
	\label{fig:isospectral-domain}
\end{figure}
\begin{figure}
	\centering
	\begin{minipage}[t]{0.49\textwidth}
		\centering
		\includegraphics[scale=0.8]{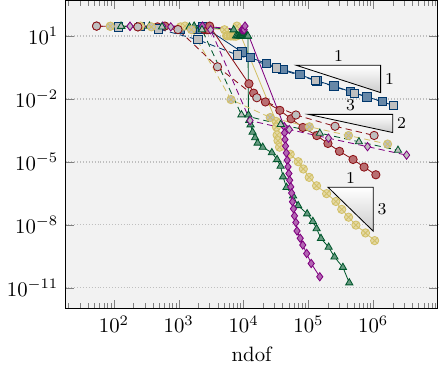}
	\end{minipage}\hfill
	\begin{minipage}[t]{0.49\textwidth}
		\centering
		\includegraphics[scale=0.8]{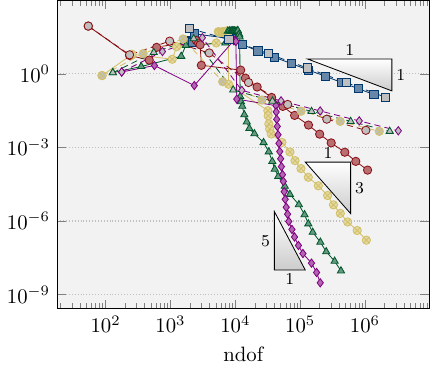}
	\end{minipage}\hfill
	\caption{Convergence history plot of $\lambda(25) - \mathrm{GLB}(25)$ (left) and $\eta^2$ (right) for polynomial
	degrees $p$ from \Cref{fig:legend} on uniform (dashed line) and adaptive (solid line) triangulations of the
isospectral domain in \Cref{sec:isospectral-domain}.}
	\label{fig:isospectral-domain-index-25}
\end{figure}
\begin{figure}
	\centering
	\begin{minipage}[t]{0.49\textwidth}
		\centering
		\includegraphics[scale=0.8]{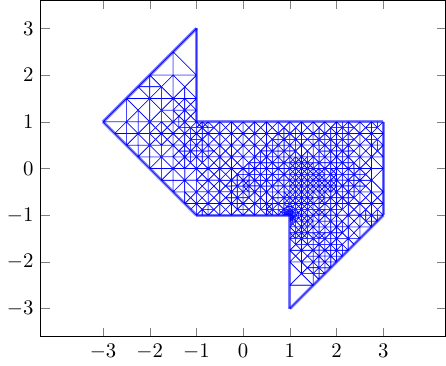}
	\end{minipage}\hfill
	\begin{minipage}[t]{0.49\textwidth}
		\centering
		\includegraphics[scale=0.8]{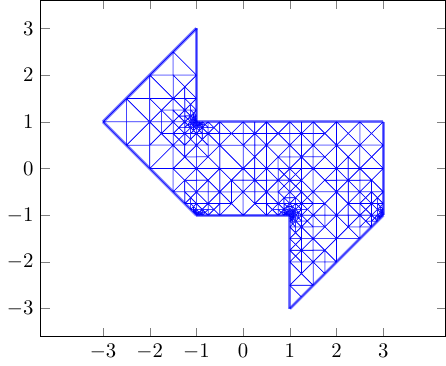}
	\end{minipage}\hfill
	\caption{Adaptive triangulations of the isospectral domain in \Cref{sec:isospectral-domain} into 1342 triangles for
		$p = 0$ (left) and into 1311 triangles for $p = 3$ (right) for the approximation to $\lambda(1)$.}
	\label{fig:isospectral-drums-triangulation}
\end{figure}
The isospectral drums are pairs of non-isometric domains with identical spectrum of the Laplace operator. This
subsection considers the domain $\Omega$ shown in \Cref{fig:initial-triangulation}.b from \cite{GordonWebbWolpert1992};
the reference values $\lambda(1) = 2.53794399980$ and $\lambda(25) = 29.5697729132$ are from \cite{BetckeTrefethen2005}
and \cite{Driscoll1997}. \Cref{fig:isospectral-domain} shows the suboptimal convergence rate $2/3$ for $\lambda(1) -
\mathrm{GLB}(1)$ and $\eta^2$ for the approximation of the principle eigenvalue $\lambda(1)$ on uniformly refined triangulations. The adaptive mesh-refining algorithm refines towards four singular corners (for $p = 3$) as depicted in \Cref{fig:isospectral-drums-triangulation} and recovers the optimal convergence rates $p+1$ for $\lambda(1) - \mathrm{GLB}(1)$ and $\eta^2$.
\Cref{fig:isospectral-domain-index-25} displays the empirical convergence rate $1$ for both $\lambda(25) - \mathrm{GLB}(25)$ and $\eta^2$ in case $p = 0$, while it is the expected rate $2/3$ for $p\ge 1$ in the presence of a typical corner singularity in the eigenfunction. We conjecture that the singular contribution to the corresponding eigenfunction in this particular example has a very small coefficient and the reduced asymptotic convergence rate $2/3$ is therefore barely visible unless a very high accuracy is reached (e.g., absolute error in the eigenvalues much smaller than $5\times 10^{-4}$). 
The adaptive mesh-refining algorithm refines towards four re-entering corners and recovers the optimal convergence rates $p+1$ for $\lambda(25) - \mathrm{GLB}(25)$ and $\eta^2$. There are two reasons for the plateau observed in the convergence history plot of $\lambda(25)-\mathrm{GLB}(25)$ in \Cref{fig:isospectral-domain-index-25}.a. First, a larger pre-asymptotic range is expected and observed for the approximation of larger eigenvalues. Second, the condition $h_{\max}^2 \lambda_h \leq \alpha$ is not satisfied for the first triangulations, whence $\mathrm{GLB}$ is set to zero. An asymptotic behaviour is observed beyond 30,000 degrees of freedom for all displayed polynomial degrees $p = 0, \dots, 4$.
\subsection{Dumbbell-slit domain}\label{sec:dumbbellslit}
\begin{figure}
	\centering
	\begin{minipage}[t]{0.49\textwidth}
		\centering
		\includegraphics[scale=0.8]{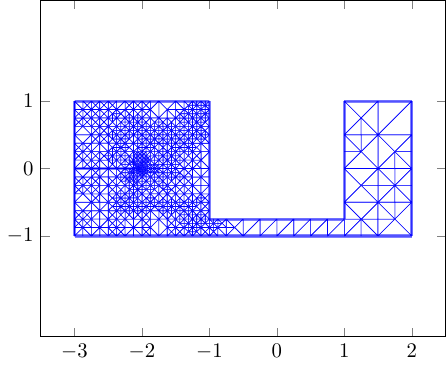}
	\end{minipage}\hfill
	\begin{minipage}[t]{0.49\textwidth}
		\centering
		\includegraphics[scale=0.8]{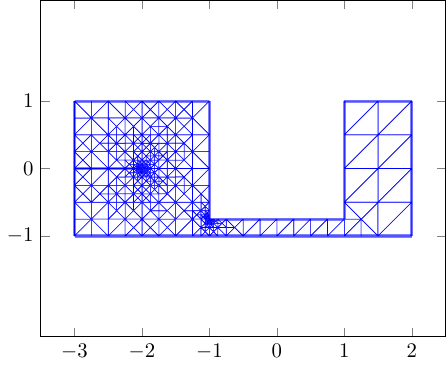}
	\end{minipage}\hfill
	\caption{Adaptive triangulations of the dumbbell-slit domain in \Cref{sec:dumbbellslit} into 1588 triangles for $p =
	0$ (left) and into 1505 triangles for $p = 3$ (right) for the approximation to $\lambda(1)$.}
	\label{fig:dumbbellslit-triangulation}
\end{figure}

\begin{figure}
	\centering
	\begin{minipage}[t]{0.49\textwidth}
		\centering
		\includegraphics[scale=0.8]{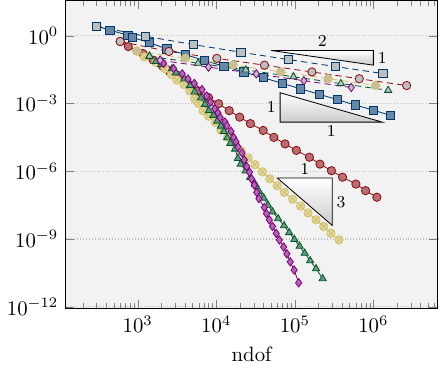}
	\end{minipage}\hfill
	\begin{minipage}[t]{0.49\textwidth}
		\centering
		\includegraphics[scale=0.8]{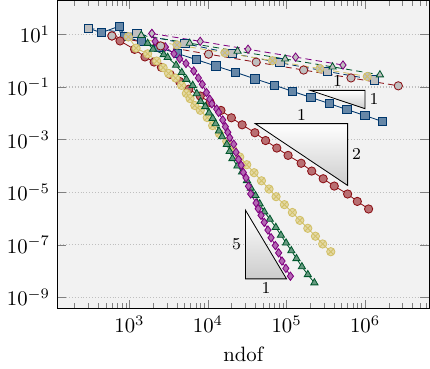}
	\end{minipage}\hfill
	\caption{Convergence history plot of $\lambda(1) - \mathrm{GLB}(1)$ (left) and $\eta^2$ (right) for polynomial
	degrees $p$ from \Cref{fig:legend} on uniform (dashed line) and adaptive (solid line) triangulations of the
dumbbell-slit domain in \Cref{sec:dumbbellslit}.}
	\label{fig:dumbbellslit}
\end{figure}

The final example approximates the principle Dirichlet eigenvalue $\lambda(1)$ on the domain $\Omega \coloneqq (-3,2)
\times (-1,1) \setminus ((-3,-2]\times\{0\} \cup [-1,1] \times [-3/4,1))$ displayed in
\Cref{fig:initial-triangulation}.c. This is a modification of the numerical example in \cite[Section 4.2]{CP21b}. The
reference value $\lambda(1) = 8.367702430882$ stems from an adaptive computation with the polynomial degree $p = 5$. The
adaptive algorithm refines towards the reentrant corners at $(-1,-3/4)$ and $(-2,0)$ as displayed in
\Cref{fig:dumbbellslit-triangulation}, while the triangles in the subdomain $(1,2) \times (-1,1)$ remain unchanged for $p \geq 1$.
Hence, there may be no reduction of the maximal mesh-size $h_{\max}$.
\Cref{fig:dumbbellslit} displays suboptimal convergence rate $0.5$ for the errors $\lambda(1) - \mathrm{GLB}(1)$ and
$\eta^2$ for $p = 0,\dots,4$. The adaptive
mesh-refining recovers the optimal convergence rates $p+1$.

\subsection{Conclusions}
The computer experiments provide empirical evidence for optimal convergence rates of the adaptive mesh-refining algorithm.
The ad hoc choice $\alpha = 1/2$ is robust in all computer experiments.
For $\beta = \alpha/\sigma_2^2$, the computable condition $\sigma_2^2 h_{\max}^2 \lambda_h(j) \leq \alpha$ leads to confirmed lower eigenvalue bounds and holds on triangulations into right-isosceles triangles, whenever the maximal mesh-size $h_{\max}$ satisfies $
\lambda_hh_{\max}^2 \leq \alpha \pi^2$.
In all displayed numerical benchmarks, $\lambda_h$ is a lower eigenvalue bound of $\lambda$ even for $\lambda_hh_{\max}^2 > \alpha \pi^2$.
The computed (but otherwise undisplayed) efficiency indices $7 {\times} 10^{-2} \le I \coloneqq |\lambda - \lambda_h|\eta^{-2} \leq 4 {\times} 10^{-3}$ range in the numerical examples
from $.07$ to $4 {\times} 10^{-3}$ for an asymptotic range $2 {\times} 10^4\le \mathrm{ndof} \le 10^5$; the quotient $I$ decreases for larger polynomial degree $p$.
The overall numerical experience provides convincing evidence for the efficiency and reliability of the stabilization-free a posteriori error estimates of this paper.
Higher polynomial degrees $p$ lead to significantly more accurate lower bounds and clearly outperform lowest-order discretizations.
\printbibliography
\appendix
\renewcommand{\thetable}{S\arabic{table}}
\setcounter{section}{0}
\renewcommand{\thesection}{\arabic{section}}
\setcounter{equation}{0}
\setcounter{figure}{0}
\setcounter{table}{0}
\makeatletter
\renewcommand{\theequation}{\arabic{equation}}
\renewcommand{\thefigure}{S\arabic{section}.S\arabic{figure}}
\newpage
\section*{Appendix: On $p$-robustness of constants \\in refined $H^1$ stability estimates}

This appendix provides details of the proof of Theorem 2.1 in the paper with focus on the 
constants $C_{\mathrm{st,1}}, C_\mathrm{st,2}$ and their dependence on the polynomial degree $p \in \mathbb{N}_0$ in three space dimensions. 
\subsection*{Overview}%
\label{sub:Overview}
Let $\vertiii{\bullet} \coloneqq \|\nabla \bullet\|_{L^2(T)}$ abbreviate the seminorm in the Sobolev space $H^1(T)\coloneqq H^1(\mathrm{int}(T))$
and let $\Pi_p$ denote the $L^2$ projection onto the space $P_p(T)$ of polynomials of total degree at most $p \in
\mathbb{N}_0$ for a
fixed tetrahedron $T \subset \R^3$. 
For any Sobolev function $f \in H^1(T)$, the Galerkin projection $G f \in P_{p+1}(T)$ is the unique polynomial of degree at most $p+1$ with the prescribed integral mean $\Pi_0 G f = \Pi_0 f$ and the orthogonality $\nabla (f - Gf) \perp \nabla P_{p+1}(T)$ in $L^2(T;\R^3)$. The constants $C_\mathrm{st,1}$ and $C_\mathrm{st,2}$ are the best possible constants in the stability estimates
\begin{align}
	\vertiii{(1 - \Pi_{p+1}) f} &\leq C_\textrm{st,1}\|(1 - \Pi_p)\nabla f\|_{L^2(T)}\quad\text{for all }f\in H^1(T),
	\label{ineq:stability-1-appendix}\\
	\vertiii{(1 - G) f} &\leq C_\textrm{st,2}\|(1 - \Pi_p)\nabla f\|_{L^2(T)}\quad\text{for all }f\in H^1(T).
	\label{ineq:stability-2-appendix}
\end{align}
Theorem 2.1 asserts the following properties of $C_\mathrm{st,1}$ and $C_\mathrm{st,2}$.
\begin{enumerate}[wide]
	\item[(A)] There exist positive constants $1 \leq C_\mathrm{st,2} \leq C_\mathrm{st,1} < \infty$ that satisfy
	\eqref{ineq:stability-1-appendix}--\eqref{ineq:stability-2-appendix}. The constants $C_\mathrm{st,1}$ and $C_\mathrm{st,2}$ are independent of the diameter $h_T$ of $T$.
	\item[(B)] $C_\mathrm{st,2}$ is $p$ robust, i.e., $C_\mathrm{st,2}$ is uniformly bounded for all $p \in \mathbb{N}_0$.
	\item[(C)] $C_\mathrm{st,1} \approx \sqrt{p+1}$ is not $p$ robust.
\end{enumerate}
The proof of $\|\Pi_{p+1}\|\leq C_{\rm st, 1}\leq C_{\rm st, 2} \|\Pi_{p+1}\|\lesssim \sqrt{p+1}$ is already explained
in the paper and $\sqrt{p+1}\lesssim \|\Pi_{p+1}\|$ 
is established below in C.
\subsection*{A.\ Proof of existence}%
\label{sub:Proof of (A)}
The two assumptions (H1)--(H2) from \cite[Theorem 3.1]{CarstensenZhaiZhang2020} imply the existence of the constant $C_\mathrm{st,1} < \infty$ in \cite[Theorem 2.3]{CarstensenErnPuttkammer2021}. The $L^2$ orthogonality $\nabla (1 - G) f \perp \nabla P_{p+1}(T)$ implies $\vertiii{(1 - G) f} \leq \|\nabla (1 - \Pi_{p+1}) f\|_{L^2(T)}$ for all $f \in H^1(T)$, whence $C_{\mathrm{st,1}} \leq C_\mathrm{st,2} < \infty$.
The best approximation property of the $L^2$ projection $\Pi_p$ proves $\|(1 - \Pi_p) \nabla f\|_{L^2(T)} \leq
\vertiii{(1 - G) f}$ and, therefore, $1 \leq C_\mathrm{st,1}$.
Notice that (A) holds in any space dimension.\qed

\subsection*{B.\ Proof of $p$ robustness of $C_{\rm st, 2}$}%
\label{sub:Proof of B}
Let $\mathrm{N}_p(T) \coloneqq P_p(T;\R^3) \oplus (P_p^\mathrm{hom}(T;\R^3) \times x) = P_p(T;\R^3) \oplus \{q \in P_{p+1}^{\mathrm{hom}}(T;\R^3): x \cdot q(x) = 0 \text{ for all } x \in T\}$ denote the first-kind N\'ed\'elec finite element space with the space $P_p^\mathrm{hom}(T;\R^3)$ of homogenous polynomials of (total) degree $p$.
Since $P_p(T;\R^3) \subset \mathrm{N}_p(T)$, the $L^2$ projection $\Pi_\mathrm{N}$ onto $\mathrm{N}_p(T)$ satisfies $\|(1 - \Pi_\mathrm{N}) \nabla f\|_{L^2(T)} \leq \|(1 - \Pi_p) \nabla f\|_{L^2(T)}$ for all $f \in H^1(T)$. Hence, the existence of a constant $C(T)$ independent of $p$ and $\diam(T)$ such that
\begin{align}
		\vertiii{f - Gf} \leq C(T)\|(1 - \Pi_\mathrm{N})\nabla f\|_{L^2(T)} \quad\text{for all } f \in H^1(T)\label{ineq:stability-C(T)}
\end{align}
implies (B). 	Given any $f \in H^1(T)$, abbreviate $q_\mathrm{N} \coloneqq \Pi_\mathrm{N} \nabla f \in \mathrm{N}_p(T)$ and
observe $r_\RT
\coloneqq \mathrm{curl}\,q_\mathrm{N} \in \RT_p(T)$ with $\div\,r_\RT = 0$, e.g., from \cite[Lemma
15.10]{ErnGuermondII2021}, \cite[Eq.~(2.3.62)]{BoffiBrezziFortin2013}, or \cite[Lemma 5.40]{Monk2003}.
It goes back to \cite{CostabelMcIntosh2010} to define a Bogovski\v{\i}-type integral operator as a pseudo-differential operator of order $-1$ of a H\"ormander class $S^{-1}_{1,0}(\R^n)$ that leads to right-inverses for differential operators.
In particular, there exists a bounded linear operator $R^{\mathrm{curl}} : H^{-1}(T;\R^3) \to L^2(T;\R^3)$ such that $R_\mathrm{N} \coloneqq R^{\mathrm{curl}} r_\RT \in \mathrm{N}_p(T)$ satisfies $\mathrm{curl}\, R_\mathrm{N} = r_\RT$. Since $R_\mathrm{N} - q_\mathrm{N} \in \mathrm{N}_p(T)$ is curl-free by design, $R_\mathrm{N} - q_\mathrm{N} = \nabla \psi$ is the gradient of some function $\psi \in H^1(T)$ in the tetrahedron $T$.
The structure of $\mathrm{N}_p(T)$ enforces $\psi \in P_{p+1}(T)$ (cf.~\cite[Lemma 15.10]{ErnGuermondII2021} and \cite[Lemma 5.28]{Monk2003} for the proof).
Recall that $\nabla G f$ is the best-approximation of $\nabla f$ in $\nabla P_{p+1}(T)$ and deduce (from $\nabla P_{p+1}(T) \subset P_p(T;\R^3) \subset N_p(T)$) that it is also the best-approximation of $q_\mathrm{N} = \Pi_\mathrm{N} \nabla f$. Hence,
\begin{align}
	\|\Pi_\mathrm{N} \nabla f - \nabla Gf\|_{L^2(T)} \leq \|q_\mathrm{N}\ + \nabla \psi\|_{L^2(T)} = \|R_\mathrm{N}\|_{L^2(T)}.
	\label{ineq:RN-bound}
\end{align}
The operator norm $\|R^{\mathrm{curl}}\|$ of $R^{\mathrm{curl}}$ allows for $\|R_\mathrm{N}\|_{L^2(T)} \leq
\|R^\mathrm{curl}\| \vertiii{r_\RT}_*$ with the norm $\vertiii{\bullet}_*$ in the dual space $H^{-1}(T;\R^3)$ of $H^1_0(T;\R^3)$ (endowed with the $H^1$ seminorm $\vertiii{\bullet}$), i.e.,
\begin{align*}
	\vertiii{r_\RT}_* = \sup_{v \in H^1_0(T;\R^3)\setminus\{0\}} \int_T r_\RT \cdot v \d{x}/\vertiii{v}.
\end{align*}
Recall $r_\RT = \mathrm{curl}\,q_\mathrm{N}$. An integration by parts and $\mathrm{curl}\,\nabla f = 0 \in L^2(T)$ {provide}
\begin{align*}
	\int_T r_\RT \cdot v \d{x} = \int_T \mathrm{curl}(q_\mathrm{N} - \nabla f) \cdot v \d{x} = \int_T (1 - \Pi_\mathrm{N}) \nabla f \cdot \mathrm{Curl}\,v \d{x}
\end{align*}
for any $v \in H^1_0(T;\R^3)$. This, a Cauchy inequality, and the estimate $\|\mathrm{Curl}\,v\|_{L^2(T)} \leq
2\vertiii{v}/\sqrt{3}$ {reveal} $\vertiii{r_\RT}_* \leq 2\|(1 - \Pi_\mathrm{N})\nabla f\|_{L^2(T)}/\sqrt{3}$. Hence
\eqref{ineq:RN-bound} implies \begin{align*}\|\Pi_\mathrm{N} \nabla f - \nabla Gf\|_{L^2(T)} \leq
	2\|R^{\mathrm{curl}}\|\|(1 - \Pi_\mathrm{N})\nabla f\|_{L^2(T)}/\sqrt{3}.\end{align*} This and the Pythagoras
theorem result in
\begin{align*}
	\vertiii{(1 - G)f}^2 &= \|(1 - \Pi_\mathrm{N})\nabla f\|_{L^2(T)}^2 + \|\Pi_\mathrm{N} \nabla f - \nabla G f\|_{L^2(T)}^2\\
	&\leq (1 + 4\|R^{\mathrm{curl}}\|^2/3) \|(1 - \Pi_\mathrm{N})\nabla f\|_{L^2(T)}^2.
\end{align*}
This proves \eqref{ineq:stability-C(T)} with $C(T) \coloneqq \sqrt{1 + 4\|R^{\mathrm{curl}}\|^2/3}$ and, therefore,
(B). More details on $\|R^{\mathrm{curl}}\|$ and further applications can be found in \cite[Section
3]{CostabelMcIntosh2010}, \cite[Section 2]{Hiptmair2009}, and \cite[Lemma 6.4]{MelenkRojik2020}.

{
An alternative proof of \eqref{ineq:stability-C(T)} involves the main result of \cite{ChaumontErnVohralik2020} and was kindly provided by A.~Ern in private communications from 03/08/2022.
For $v \coloneqq \nabla f \in H(\mathrm{curl}, T)$ with $\mathrm{curl}\, v = 0$, let
$v_h^*$ (resp.~$w_h^*$) denote the minimizer of $\|v - v_h\|_{L^2(\Omega)}$ among $v_h \in \mathcal{K} \coloneqq \{v_h \in N_p(T) : \mathrm{curl}\,v = 0\}$ (resp.~$\|v - w_h\|_{L^2(\Omega)}$ among $w_h \in N_p(T)$).
The $L^2$ orthogonality $w_h^* - v_h^* \perp \mathcal{K}$ from the Euler-Lagrange equations associated with these minimization problems
implies that the difference
$w_h^* - v_h^*$ minimizes the functional $\|z_h\|_{L^2(T)}$ among all $z_h \in w_h^* + \mathcal{K}$.
Invoking the results of \cite{CostabelMcIntosh2010}, it is known from \cite[Theorem 2]{ChaumontErnVohralik2020} that
\begin{align*}
	\|w_h^* - v_h^*\|_{L^2(T)} = \inf_{\substack{z_h \in \mathrm{N}_p(T)\\\mathrm{curl}\,z_h = \mathrm{curl}\, w_h^*}}\|z_h\|_{L^2(T)}
	\leq \widetilde{C}(T)\inf_{\substack{z \in H(\mathrm{curl},T)\\\mathrm{curl}\, z = \mathrm{curl}\, w_h^*}}\|z\|_{L^2(T)}
\end{align*}
with a $p$-robust constant $\widetilde{C}(T) > 0$.
Since $\mathrm{curl}(w_h^* - v) = \mathrm{curl} w_h^*$, we infer $\|w_h^* - v_h^*\|_{L^2(T)} \leq \widetilde{C}(T)\|v - w_h^*\|_{L^2(T)}$. This and a triangle inequality imply
\begin{align*}
	\|v - v_h^*\|_{L^2(T)} \leq (1 + \widetilde{C}(T))\|v - w_h^*\|_{L^2(T)}.
\end{align*}
This is \eqref{ineq:stability-C(T)} with $C(T) \coloneqq 1 + \widetilde{C}(T)$ because $w_h^* = \Pi_\mathrm{N} \nabla f$ by design and $v_h^* = \nabla G f$ from $\mathcal{K} = \nabla P_{p+1}(T)$.\qed
}
	\subsection*{C. Lower growth $\sqrt{p+1}\lesssim \|\Pi_{p+1}\|$}%
	\label{sec:Supplement C: Proof of the lower growth {$sqrt{p}lesssim |Pi_{p}$}}
	While a compactness argument in \cite[Theorem 3.1]{CarstensenZhaiZhang2020} leads to the existence of
	$C_\mathrm{st,1}$, the dependence of $C_\mathrm{st,1}$ on the polynomial degree $p$ remained obscured and only an upper bound for $p = 0$ was given.
	The proof of Theorem 2.1 in the paper establishes 
	\begin{align*}
		C_\mathrm{st,1}\approx \|\Pi_{p+1}\|\coloneqq \sup_{\phi\in
			H^1(T)\setminus\R}\frac{\vertiii{\Pi_{p+1}\phi}}{\vertiii{\phi}}.
	\end{align*}
	An upper bound $\|\Pi_{p+1}\|\lesssim \sqrt{p+1}$ of the growth of the $H^1$ stability constant of the $L^2$ projection
	is known
	from \cite[Sec.\ 5]{melenk_stability_2013} and \cite{wurzer2010stability}.
	The remaining parts of this appendix therefore consider the reverse direction $\sqrt{p+1}\lesssim \|\Pi_{p+1}\|$ for a
	tetrahedron and
	depart with a motivating classical result in 1D.
	For simplicity, the following presentation applies an index shift and discusses $\|\Pi_p\|\approx \sqrt{p}$ for arbitrary
	$p\geq1$.
	\subsubsection*{Lower bound in 1D}%
	\label{sub:1D case}
	In one space dimension, $\|\Pi_p\|\approx\sqrt{p}$ is established, e.g., in \cite[Theorem
	2.4]{canuto_approximation_1982} and \cite[Remark 3.5]{bernardi_polynomial_1992}.
	Let $L_k$ for $k \in \mathbb{N}_0$ denote the Legendre polynomials in the reference interval $I\coloneqq(-1,1)$.
	Then $L_k$ satisfies, for all $k \in \mathbb{N}_0$,
	\begin{align}
		\widehat L_k(x) &\coloneqq \int_{-1}^x L_k(t) \d{t} = \frac{L_{k+1}(x)-L_{k-1}(x)}{2k+1},\label{eqn:B1}\\
		\|L_k\|_{L^2(I)}^2 &= \frac{2}{2k+1}\leq \frac{1}{k}, \quad\text{and}\quad\|\nabla L_k\|_{L^2(I)}^2 = k(k+1)\label{eqn:B2}
	\end{align}
	with the convention $L_{-1} \equiv 0$ in $I$,
	cf., e.g., \cite[Eqns (3.11), (3.12), (5.3)]{bernardi_spectral_1997}.
	The pairwise $L^2$ orthogonality of $L_k$ and \eqref{eqn:B1}--\eqref{eqn:B2} lead to
	\begin{align*}
		\|\nabla \Pi_p \widehat L_p\|_{L^2(I)}^2 &= \frac{\|\nabla L_{p-1}\|_{L^2(I)}^2}{(2p+1)^2} =
		\frac{p(p-1)}{(2p+1)^2}\approx 1,\\
		\|\nabla \widehat L_{p}\|_{L^2(I)}^2 &= \|L_{p}\|_{L^2(I)}^2 = \frac{2}{2p+1}\approx p^{-1},
	\end{align*}
	whence
	$\|\nabla \Pi_{p} \widehat L_{p}\|_{L^2(I)}\approx \sqrt{p}\|\nabla \widehat L_{p}\|_{L^2(I)}$ for all $p \geq 1$.
	This proves $\sqrt{p} \lesssim \|\Pi_p\|$ in 1D.
	
	A similar result holds for the $L^2$ projection $\widetilde \Pi_{p}:L^2(D)\to Q_{p}(D) \cong P_{p}(I)^n$ onto the space of tensor-product
	polynomials on the $n$-cube $D\coloneqq I^n=(-1,+1)^n$ \cite{canuto_approximation_1982}. 
	For simplicity and because the arguments carry over to triangles as well, the following proof
	considers simplices in $n=3$ dimensions only.\medskip
	
	\noindent\emph{Proof of $\sqrt p\lesssim \|\Pi_p\|$ for $n=3$}.
	Let $p\in\mathbb N$ be arbitrary and let $F:Q\to T$ denote the coordinate transformation 
	\begin{align*}
		F(\eta_1, \eta, \eta_3) \coloneqq \left(\frac{(1+\eta_1)(1-\eta_2)(1-\eta_3)}{4}-1,
		\frac{(1+\eta_2)(1-\eta_3)}{2}-1, \eta_3\right)
	\end{align*}
	from the cube $Q\coloneqq(-1,1)^3$ onto the reference tetrahedron $T \coloneqq \mathrm{conv}\{(-1,-1,-1),\allowbreak (1,\allowbreak -1,\allowbreak -1), (-1, 1, -1), (-1, -1, 1)\}$ with the Jacobian $J_F$ and $\mathrm{det} J_F=(1-\eta_2)(1-\eta_3)^2/8$, see, e.g.,
	\cite{sherwin_karniadakis_orthogonal_basis_3D} and \cite[Section 3]{melenk_stability_2013} for a derivation. An integration by substitution leads, for any $f \in L^1(Q)$, to
	\begin{align}\label{eqn:int_by_subs}
		\int_{T} f\circ F^{-1} \d{x} = \frac{1}{8}\int_Q (1-\eta_2)(1-\eta_3)^2 f \d{(\eta_1, \eta_2,\eta_3)}.
	\end{align}
	Define $\varphi(\eta_2,\eta_3)\coloneqq(1-\eta_2)(1-\eta_3)$, $U_p(\eta_1, \eta_2,\eta_3)\coloneqq
	\varphi(\eta_2,\eta_3)^{p-1}\widehat L_p(\eta_1)$, and $\widetilde U_p\coloneqq U_p\circ F^{-1}\in L^2(T)$ for $p \geq 1$.
	The chain rule $\nabla\widetilde U_p = J_F^{-\top}\nabla U_p\circ F^{-1}$ for the gradient and \eqref{eqn:B1}
	provides $(\nabla \widetilde U_p)\circ F = \varphi(\eta_2, \eta_3)^{p-2} M(\eta_1, \eta_2) G(\eta_1)$ with
	\begin{align*}
		M(\eta_1, \eta_2)\coloneqq \begin{pmatrix}
			4&0&0\\
			2(1+\eta_1)&2&0\\
			2(1+\eta_1)&(1+\eta_2)&(1-\eta_2)
		\end{pmatrix}&&\text{and}&& G(\eta_1)\coloneqq \begin{pmatrix}
			L_p(\eta_1)\\(p-1)\widehat L_p(\eta_1)\\(p-1)\widehat L_p(\eta_1)
		\end{pmatrix}.
	\end{align*}
	A Cauchy inequality in $\R^3$ proves
	\begin{align*}
		|\nabla\widetilde U_p|^2\circ F &= \varphi(\eta_2, \eta_3)^{2p-4}\sum^{3}_{j=1} \Big(\sum^{3}_{k=1}
		M_{jk}G_k\Big)^2
		\leq 3\varphi(\eta_2, \eta_3)^{2p-4}\sum^{3}_{k=1} \Big(\sum^{3}_{j=1} M_{jk}^2\Big)G_k^2.
	\end{align*}
	This, the integration by substitution formula \eqref{eqn:int_by_subs},
	$R_p\coloneqq\int_{-1}^1\int_{-1}^1(1-\eta_2)^{2p-3}(1-\eta_3)^{2p-2}\d{(\eta_2,\eta_3)}\in \mathbb R$, and $|\eta_j|\leq 1$ for $j=1,2,3$ and $(\eta_1,\eta_2, \eta_3)\in Q$
	show
	\begin{align}\notag
		\frac{1}{3}\vertiii{\widetilde U_p}^2 \leq&\; \frac{1}{8}\int_Q (1-\eta_2)(1-\eta_3)^2 c(\eta_2,
		\eta_3)^{2p-4}\\&\mkern-48mu\times\big((16+8(1-\eta_1)^2)L_p(\eta_1)^2 + 2(3+\eta_2^2)(p-1)^2\widehat L_p(\eta_1)^2\big)\d{(\eta_1,
			\eta_2,\eta_3)}\notag\\\notag
		\leq&\; R_p\int_{-1}^16L_p(\eta_1)^2 + (p-1)^2\widehat L_p(\eta_1)^2\d{\eta_1}\\
		=&\;R_p\big(6\|L_p\|_{L^2(I)}^2 + \|(p-1)\widehat L_p\|_{L^2(I)}^2\big).\label{eqn:B_nabla_bound}
	\end{align}
	The pairwise $L^2$ orthogonality of Legendre polynomials and \eqref{eqn:B1}--\eqref{eqn:B2} verify
	\begin{align*}
		p\|(p-1)\widehat L_p\|_{L^2(I)}^2 &= \frac{p(p-1)^2}{(2p+1)^2}\left(\|L_{p+1}\|_{L^2(I)}^2 +
		\|L_{p-1}\|_{L^2(I)}^2\right)\\
		&=\frac{p(p-1)^2}{(2p+1)^2}\left(\frac{2}{2p+3}+\frac{2}{2p-1}\right)\\
		&=\frac{(p-1)^2p}{2(p-1/2)(p+1/2)(p+3/2)}\leq
		\frac{1}{2}
	\end{align*}
	for $p\geq 1$.
	This, \eqref{eqn:B2}, and \eqref{eqn:B_nabla_bound} provide the bound $2p\vertiii{\widetilde U_p}^2\leq
	39 R_p$ for $p\geq 1$.
	It remains to control $\nabla \Pi_p \tilde U_p$ from below.
	Recall from \cite{sherwin_karniadakis_orthogonal_basis_3D} that the polynomials $\widetilde \psi_{j,k,\ell}\coloneqq \psi_{j,k,\ell}\circ F^{-1}\in
	P_{j+k+\ell}(T)$ for $j,k,\ell\in\mathbb N_0$ with 
	\begin{align*}
		\psi_{j,k,\ell}(\eta_1, \eta_2,\eta_3)\coloneqq
		L_j(\eta_1)\left(1-\eta_2\right)^j\Pkab k{2j+1}(\eta_2)\left(1-\eta_3\right)^{j+k}\Pkab\ell{2j+2k+2}
		(\eta_3)
	\end{align*}
	are $L^2(T)$ orthogonal and that $(\widetilde \psi_{j,k,\ell}\ |\ 0\leq j+k+\ell\leq p)$ forms a basis of $P_p(T)$.
	The pairwise orthogonality of the Legendre polynomials, \eqref{eqn:B1}, and \eqref{eqn:int_by_subs} imply that
	\begin{align*}%
		\left(\widetilde U_p + \frac{\widetilde\psi_{p-1,0,0}}{2p+1},
		\widetilde\psi_{j,k,\ell}\right)_{L^2(T)}=\left((1-\eta_2)^{p}(1-\eta_3)^{p+1}	
		\frac{L_{p+1}(\eta_1)}{2p+1}, \psi_{j,k,\ell}\right)_{L^2(Q)}=0
	\end{align*}
	vanishes
	for all $k,\ell\in\mathbb N$ and $j\leq p$. Consequently, $$(2p+1)\Pi_p\widetilde U_p =-\widetilde\psi_{p-1,0,0}\in P_{p-1}(T).$$
	This, the chain rule for partial derivatives, and \eqref{eqn:B2}--\eqref{eqn:int_by_subs} show
	\begin{align}
		\left\|\frac{\partial}{\partial x}\Pi_p \widetilde U_p\right\|_{L^2(T)}^2 &=
		\frac{4^2R_p}{8(2p+1)^2}\int_{-1}^1\left(\frac{\mathrm{d}\phantom{A}}{\!\d{\eta_1}}L_{p-1}(\eta_1)\right)^2\!\!\!\d{\eta_1}\nonumber\\
		&=\frac{2R_p}{(2p+1)^2}\|\nabla
		L_{p-1}\|_{L^2(I)}^2=\frac{2R_p}{(2p+1)^2}p(p-1).
		\label{ineq:proof-h1-stability-L2-projection-upper-bound}
	\end{align}
	The term $2p(p-1)(2p+1)^{-2}\geq0$ is monotonically increasing in $p\geq 1$ and
	bounded from below by $4/25$ for $p\geq 2$.
	Thus, \eqref{ineq:proof-h1-stability-L2-projection-upper-bound} and $2p\vertiii{\widetilde U_p}^2\leq
	39 R_p$ provide 
	\begin{align*}
		\frac{8}{975}p\vertiii{\widetilde
			U_p}^2\leq\frac{4}{25}R_p\leq\left\Vert \frac{\partial}{\partial x}\Pi_p\widetilde
		U_p\right\Vert_{L^2(T)}^{2} \leq\vertiii{\Pi_p\widetilde U_p}^2
	\end{align*}
	for all $p\geq 2$, whence $\sqrt{p} \lesssim \|\Pi_p\|$ on the reference tetrahedron $T$.
	This and a scaling argument with an affine transformation concludes the
	proof for a general tetrahedron.\qed
	
	\subsection*{Acknowledgements}
	The authors gratefully thank Prof.\ Markus Melenk (Vienna University of Technology) for the discussion about the $H^1$ stability of the $L^2$
	projection that eventually led to the proof of (C)
	{and Prof.\ Alexandre Ern (CERMICS, ENPC) for his alternative proof of the $p$-robustness of $C_{\mathrm{st,2}}$ in (B).}
\end{document}